\documentclass[12pt]{article}

\usepackage{amsmath}
\usepackage{latexsym}
\usepackage{graphicx}
\usepackage{epsf}

\usepackage{amsmath}
\usepackage{amssymb}
\usepackage{exscale,cmmib57}

\chardef\bslash=`\\ % p. 424, TeXbook
\makeatletter
\def\verbatim{\interlinepenalty\@M \@verbatim
  \leftskip\@totalleftmargin\advance\leftskip2pc
  \frenchspacing\@vobeyspaces \@xverbatim}
\makeatother
\newtheorem{thm}{Theorem}[section]
\newtheorem{cor}[thm]{Corollary}
\newtheorem{lem}[thm]{Lemma}
\newtheorem{prop}[thm]{Proposition}

\newtheorem{defn}[thm]{Definition}

 \newtheorem{exmp}{Example}[section]

\newtheorem{remark}{Remark}[section]

\newcommand\qed{\mbox{$\Box$}}

 \newenvironment{Proof}
     {{\bf Proof:}}
     {\qed}

 \newcommand{\func}[1]{\mbox{\rm #1}\,}
 \newcommand{\funcc}[1]{\mbox{\rm #1}}
 
 \newcommand{\e}{\epsilon}
 \newcommand{\la}{\langle}
 \newcommand{\ra}{\rangle}
 \newcommand{\R}{\mathbb{R}}

 \newcommand{\C} {\mathbb{C}}
 
\makeatletter \@addtoreset{equation}{chapter} \makeatother
\begin{document}

\title{Continuous-time performance limitations for overshoot and related tracking measures }
\author{R. B. Wenczel and R. D. Hill\thanks{
Department of Mathematical and Geospatial Sciences, RMIT University, GPO Box 2476V,
Melbourne, Victoria, 3001, Australia. e-mail: robert.wenczel@rmit.edu.au, r.hill@rmit.edu.au  }}
\maketitle

\begin{abstract}

A dual formulation for the problem of determining
absolute performance limitations on overshoot, undershoot, maximum amplitude
and fluctuation minimisation for continuous-time feedback systems is
constructed. Determining, for example, the minimum possible overshoot attainable by all
possible stabilising controllers is an optimisation task that cannot be
expressed as a minimum-norm problem. It is this fact, coupled with the
continuous-time rather than discrete-time formulation, that makes these problems
challenging.
  We extend previous results to include more general reference functions, and derive new results (in continuous time)
on the influence of pole/zero locations on achievable time-domain performance.

\end{abstract}

\section{Introduction}\label{sec:cont:intro}
In this paper we study the problem of finding limits on the performance of   the error-response performance, for a specific input, for lumped continuous-time {\scriptsize SISO} systems.  We are trying to find the best possible tracking performance achievable by a lumped (or rational) {\scriptsize BIBO}--stabilising feedback controller. The theory to be presented is applicable to a large class of performance measures of practical significance, including overshoot for example. The simultaneous imposition of hard bounds on the output, in conjunction with overshoot minimisation, can also be handled, allowing consideration of the inevitable trade-off between rise-time and overshoot performance. These issues have been examined in a discrete-time setting in \cite{HEHW} etc. The investigation of time-domain performance limitations is mathematically more challenging in continuous time than in discrete time. In the continuous-time case it is harder to derive conditions under which the performance limit for rational controllers is the same as for non-rational controllers ({\em rational approximability\/}).

The usual approach to these questions, is to express this performance--limit in the form of  an optimization problem (to be computed via duality methods), over an {\em extended\/} ambient signal--space,  usually a standard {\em dual\/} Banach space, to ensure
{\em existence of optimal elements\/}
and   adequate duality properties. In general,
these optimal elements will be non-rational, and hence suffer from implementability issues. It is then of importance to find if optimal performance may be approached using rational controllers.
 Since these nonrational elements may be of doubtful practical utility, one may well abandon the search for optima, and consider only the class of rational controllers
in the optimization. This will permit   more flexibility in the choice of  the ambient Banach space since it no longer need be a dual space.

In earlier works this underlying space was chosen to be $L^\infty$ (dual of $L^1$), but this choice of space  does not, of itself, enforce the zero--steady--state--error condition.
We shall incorporate this condition by using
  the space $C_0$
(of continuous  functions of time that tend to zero at infinity) as our choice---it does the trick, but it is not a dual, so that existence of optimal elements is not assured. This choice  of space  is dictated by the continuity properties of the basic performance measures we study here.

This is the approach followed in this paper. Moreover, to side-step the rational--approximation issue, we consider
 directly   the set of error signals due to the {\em rational\/} stabilizing controllers, and
 the optimization over its closure. Under the conditions we assume in this paper, passage to this closure has no effect on the limit-of-performance.
 This   clarifies the relation between the limit-of-performance  and rational approximation, but   requires an explicit characterization of the closure of the feasible set of
error--signals. However, it is found that the closure itself has pleasant form, well-suited to  application of the duality theory.
That is, by exploiting the structure of this set we have that
 rational approximability holds by definition.

Further, we consider a general class of performance functionals and set up a general
duality framework for the analysis of such problems (in the style of \cite{HEHW}  for the case of discrete time). These functionals include as special cases the well-known  classical criteria of overshoot, undershoot and others. This treatment also allows a general fixed input (not just a step). The classical criteria just mentioned, are continuous with respect to the supremum--norm
on the signal space---this fact motivating our choice of $C_0$ as the ambient space.

  In the literature to date the optimization is performed over error signals in the bigger space $L_\infty$, with attendant consideration of rational approximability.
    In  \cite{wang-sznaier-94, wang-sznaier-96,  wang-sznaier-97}  the set of candidate controllers is expanded, moreover,  to include those that may not be {\scriptsize BIBO}--stabilising, and the limiting error is
    unconstrained.   Asymptotic stability and zero steady state error are then enforced through selection of suitable output weighting.

    Use of the space $L_\infty$ to formulate the primal has the advantage that the optimal solution for the primal problem is guaranteed to exist for the $L_\infty$ norm minimisation problem.
    For the more general cost functions considered in this paper, however, it is clearly not to be expected that an optimal solution will in general exist, neither in $L_\infty$, nor in $C_0$.

  We provide a
 derivation of a dual formulation covering a wide range of problems, extending known results. The continuous time $L_\infty$ norm minimisation problem for a fixed input was first considered in \cite{DP2},
but rational  approximability was not addressed. Miller \cite{Miller} gave a rational approximation result for response to a step reference input. The construction of rational suboptimal controllers for the continuous time $L_\infty$ norm minimisation problem has been considered in \cite{wang-sznaier-96}, \cite{wang-sznaier-94} and \cite{Halp}. Yoon \cite{MGY} and \cite{YK} extended the class of optimisation criteria, by considering convex combinations of overshoot and undershoot (with a step-input) and in these latter papers is found the first application of conjugate--duality techniques to continuous-time control.

     The constraints for the dual formulation presented here can be interpreted as arising from an open--loop dynamic system.  These structural insights are exploited to derive new results on the influence of plant pole/zero locations, or rise--time constraints on achievable performance, giving results of identical form to those obtained for discrete time in  \cite{HEHW}.
The results presented here are applicable for general reference inputs  (\cite{MGY, YK} assumes step-input) and to systems with more general exogenous inputs, for example a fixed disturbance entering at the plant output. For definiteness and simplicity our set-up is the tracking problem shown in Figure~\ref{fig:2}.

An important design consideration is the extent to which performance objectives are in conflict. It may not be possible, for example, to find a single controller that reduces both the $L_{\infty}$--norm of the error, and the negative--error (overshoot), in response to a step, to close to their fundamental limits. Furthermore, the extent to which there is conflict depends on the location of the unstable poles and non-minimum-phase zeros of the plant. The dual provides such information. Another trade-off important in practice is that between rise--time and transient performance measures such as overshoot and undershoot. The answer to this question also depends on the location of the poles and zeros of the plant, and again the interpretation of the dual as an open-loop dynamic system provides insight and new results. Also in this paper we extend known results on the minimization of the ``maximum-peak to minimum-peak" value of the error response, that is the difference between the maximum and minimum values of the error response; we term this quantity {\em fluctuation\/}. For some applications fluctuation minimization is of more concern than $L_{\infty}$--norm minimization. For some plants there then arises an unavoidable trade-off; either a small $L_{\infty}$--norm of the error response, or a low fluctuation, must be sacrificed.  To what extent they are in conflict depends on the location of the poles and zeros of the plant, and again it is the dual that provides answers.

\section{Mathematical Preliminaries}\label{sec:cont:mathprelim}
\subsection{basic notation}
We shall write $\R_+$ for the real interval $[0,+\infty)$. $\R[s]$ and $\R(s)$ denote respectively, the
spaces of real polynomials, and real--rational functions, in the complex variable $s$. The set of all proper stable
members of $\R(s)$ (i.e.\ those  with no pole in the closed complex right-half-plane nor at infinity)  is denoted by ${\mathcal S}$.
The Laplace transform of a function $f$, will be written as $\hat{f}$. We let $L^p(\R_+)$ ($1\le p \le \infty$) stand for
the space of Lebesgue $p$-integrable functions on $\R_+$. The space of continuous functions   $\varphi:\R_+\to \R$ for which
$\lim_{t\to\infty}\varphi(t)=0$, endowed with the supremum norm $\| \varphi   \|_\infty:= \sup_{t\ge 0}|\varphi(t)|$, shall be denoted by
$C_0(\R_+)$. For any $\alpha\in\R$, write  $C_{0,\alpha}(\R_+)$ for the subset consisting of those elements $\varphi\in C_0(\R_+)$ for
which $\varphi(0)=\alpha$. For a subset $A$ of a space $X$, the indicator function $\iota_A$ of $A$ is defined on $X$ by
$$ \iota_A(x):= \left\{ \begin{array}{cl}
                         0 & \text{ if } x\in A \\
                       +\infty &  \text{ if } x\notin A
                    \end{array}\right. \,.    $$
If $f:X\to \R\cup \{+\infty\}$, then $\func{dom}f$ denotes the set $\{x\in X\mid f(x)<+\infty\}$.
\subsection{Fenchel duality theorem}
\begin{defn} Let $X$ be a topological linear space.
For any $f:X\to  \R\cup\{ +\infty\}$, the (Young--Fenchel)
conjugate $f^*:X^*\to \R\cup\{ \pm\infty\}$ is defined by:
$$ \!\!\!\!\!\text{for all }x^*\in X^*,\qquad
f^*(x^*):=\sup_{x\in X}\,(\langle x,x^*\rangle-f(x))
$$
\end{defn}
It follows that $f^*$ is convex and weak${}^*$
lower--semicontinuous.

We shall require a Fenchel duality theorem in the following form (see
     \cite[Theorem 18(a) and Example $11'$]{Rock} for a more general formulation)
\begin{prop}\label{prop:BL}
Let $f:X\to  \R\cup\{ +\infty\} $ and $g:\R^n \to\R\cup\{ +\infty\}$ be
convex, with $X$  a locally-convex topological vector space.
  Let $A:X\to\R^n$ be bounded linear. Assume also that $g$ is
%polyhedral\footnote{that is, the set of points above its graph may be formed as the intersection of
%finitely-many half-spaces---this is satisfied for the indicator function of any one-point set in finite dimensions.}
    finite-valued at some point in   $A(\func{int dom}f)$. Then,
$$
\inf\{f(x)+g(Ax)\mid x\in X\}=\max\{-f^*(A^T\xi)-g^*(-\xi)\mid
\xi\in\R^n\}
$$
%If the infimum on the left is achieved for some $x_{opt}$,  and $x_{opt}^*$
% is an optimal dual element,  then
%\begin{equation} \label{align:cont}
%f(x_{opt})+f^{*}(x^*_{opt})=\langle x_{opt},x^*_{opt}\rangle .
%\end{equation}
\end{prop}
In our applications, we will take $g=\iota_{\{b\}}$, the indicator function for a singleton $\{b\}$,
and $f$ will be of the form $f_0+\iota_T$ for a finite-valued $f_0$ and a convex set $T$.

\subsubsection{duality for a space of continuous functions}
The application of the Fenchel theorem shall require a characterisation of the dual of the Banach space
of signals under consideration. In this paper, the ambient space of (error) signals will be $C_0(\R_+)$,
which  has a well-known dual space.
Indeed, the dual of $C_0(\R_+)$
 is isometrically isomorphic with the space $ \mathbf{M}(\R_+)$ of all
regular finite signed Borel measures
 $\mu$      on $\R_+$,
with the variation norm $\|\mu\|:=|\mu|(\R_+)$ (For formal
definitions of these properties, and a statement of the duality
for $C_0(\R_+)$ see, for example, \cite{Cohn}). The action of
$\mu\in\mathbf{M}(\R_+)$ as a  dual functional on $C_0(\R_+)$, is indicated by the pairing
$$ \la \mu,e\ra = \int_{\scriptsize \R_+}\!\! e\,d\mu $$
for $e\in C_0(\R_+)$. (Note also that this expression is well--defined
whenever $e$ is  bounded and (Borel-)measurable, since $\mu$ is a finite measure.)

\section{Problem Formulation}\label{sec:cont:prelim}
\subsection{primal feasible set}\label{sect:primfeasset:cont}
\begin{figure}[h]
\begin{center}
\begin{tabular}{c}
\epsfxsize= 9 cm
\epsffile{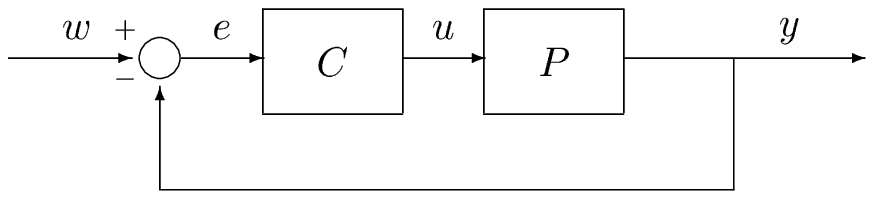}
\end{tabular}
\end{center}
\caption{\label{fig:2} A closed--loop control system
 }
\end{figure}

Consider the set ${\mathcal F}$ of all error sequences achievable with a rational stabilising controller (for the plant $P$)
in the standard one-degree-of-freedom feedback configuration, see Fig~\ref{fig:2}.   Here the plant $P$ is a scalar, linear, proper,
 finite-dimensional system. Mathematically, this condition on $P$ can be expressed by the requirement that $P\in\R(s)$ has a zero
  at $\infty$ of order $\theta_p\ge 0$ (in the sense that the degree of the denominator is the sum of
$\theta_p$ with the degree of the numerator). The reference input $w$ is assumed to be an ordinary function
 $w :\mathbf{\R_+}\to\mathbf{\R }$ (no delta-function terms) with rational Laplace transform.
  A typical linear stabilising controller is denoted $C$. The plant output is $y$, and the error signal is  $e =w-y$

  Let $\hat{n}$,
$\hat{d}$  in $\mathbf{{\mathcal S}}$  be a coprime factorization
for $P$. Then $\hat{d}(\infty):=\lim_{|s|\to\infty} \hat{d}(s)$ is
nonzero and finite, and $\hat{n}$ has  a zero of order $\theta_p$
at $\infty$, so that
$$ \lim_{|s|\to\infty}s^{\theta_p-1}\hat{n}(s)=0\text{ and }
\lim_{|s|\to\infty}s^{\theta_p}\hat{n}(s)\text{ is nonzero and
finite}\,.$$ Such a factorization for $P$ is readily found; place
$P=q/r$ for coprime polynomials $q$, $r$, with $\func{deg}
r=\func{deg}q+\theta_p$. Then
$\hat{n}:=q/(\cdot+1)^{\func{{\scriptsize deg}} r}$ and
$\hat{d}:=r /(\cdot+1)^{\func{{\scriptsize deg}} r}$  both are in
$\mathbf{{\mathcal S}}$, and by \cite[Chapter 2, fact 20]{Vidy1}, are coprime
(in $\mathbf{{\mathcal S}}$).

\smallskip

Given the reference-input $w$, then $e$ is a closed--loop error--signal (for some
 stabilizing controller $C$ for $P$)  if
and only if
$\hat{e}=\hat{w}\hat{d}(\hat{v}-\hat{q}\hat{n})=\hat{w}-\hat{w}\hat{n}(\hat{x}+\hat{q}\hat{d})$
for some $\hat{q}\in \mathbf{{\mathcal S}}$, where $\hat{x}$,
$\hat{v}$ in $\mathbf{{\mathcal S}}$ arise from the coprimeness of
$\hat{n}$ and  $\hat{d}$, and satisfy
$$ \hat{x}\hat{n}+\hat{v}\hat{d}=1 \text{ in }  \mathbf{{\mathcal S}}\,.$$

From the assumption on $w$, its Laplace transform then satisfies
$$\lim_{\func{{\scriptsize Re}}s\to +\infty}\hat{w}(s)=0\,,$$ and
$\hat{w}(s)$ has a zero of order $\theta_w\ge 1$ at infinity, implying
that $\lim_{|s|\to\infty}s^{\theta_w-1}\hat{w}(s)=0$ and that
$\lim_{|s|\to\infty}s^{\theta_w}\hat{w}(s)$ is nonzero and finite.
Hence, for $e$ as above, $\hat{e}(s)$ has a zero of order at
least $\theta_w$ at $\infty$, and from
$\hat{w}-\hat{e}=\hat{w}\hat{n}(\hat{x}+\hat{q}\hat{d})$ follows
that $\hat{e}-\hat{w}$ has a zero of order at least $\theta_p
+\theta_w$ at $\infty$.

Define our feasible set of possible error signals by the affine
set (where the star $*$ denotes convolution of functions)
$$ {\mathcal F}:= \{ e\in C_0(\R_+)\mid e=w*d*(v-q*n) \text{ for some }
q \text{ such that  } \hat{q}\in \mathbf{{\mathcal S}}\}\,.$$
 (We are constraining  $e$ to be in   $C_0(\R_+)$, as a criterion for zero steady--%
state tracking error). Then, for $\underline{\theta}:=(\theta_w,
\theta_p)$,
\begin{align*}
  {\mathcal F}\subseteq X_{\underline{\theta}} :=\{ e\in C_0(\R_+)\mid
\hat{e} &\text{ rational, and has a zero }\\
&\text{of order at least $\theta_w$ at infinity, and } \\
&\text{$\hat{e}-\hat{w}$  has zero at $\infty$ of order at least
}\theta_p+\theta_w \}.
\end{align*}
Note that $X_{\underline{\theta}}$ is an affine subspace of
$C_0(\R_+)$. By considering the partial--fraction expansions of
$\hat{e}$, it follows that any $e\in  L^1(\R_+) \cup C_0(\R_+)$  with
rational Laplace transform may be expressed in the form (for all
$t\ge 0$)

$$ e(t)=\func{Re} \sum_i c_i t^{k_i} e^{\lambda_i t} $$
where the above sum consists of a finite number of terms, and $c_i\in\C$,
$\func{Re}\lambda_i<0$, $k_i\ge 0$.

\smallskip

Introduce the notation for right half--planes
\begin{align*}
\C_+ &:= \{s\in\C\mid \func{Re}s>0\}\\
\overline{\C_+}&:= \{s\in\C\mid \func{Re}s\ge 0\}\\
  \overline{\C_+}_e&:=\overline{\C_+}\cup \{\infty\}\,.
\end{align*}
We shall now follow an analogue of the developments of
\cite{HEHW}.
 Let  $P$
have poles $p_1,\ldots, p_m$ and zeros $z_1,\ldots, z_n$ in the
right--hand plane $\overline{\C_+}$. Also, let the
reference--input $\hat{w}$ have zeros $v_1,\ldots, v_l$ in
$\overline{\C_+}$. It is assumed that none of these lie on the
imaginary  axis $\{s\mid\func{Re}s=0\}$. Further, all these
poles/zeros are assumed to be mutually distinct, and  simple.

Place 
\begin{equation}
M_{\underline{\theta}}=\left\{e\in
X_{\underline{\theta}}\,\Bigg{\vert}
\begin{array}{ll}
\hat{e}(z_i) =\hat{w}(z_i) & i=1,2,\ldots ,n \\
\hat{e}(p_i) =\,0 & i=1,2,\ldots ,m\\
\hat{e}(v_i) =\,0 & i=1,2,\ldots ,l\\
% \hat{e}(\infty)=0\text{ of order}\geq\theta_w, &   \\
% (\hat{e}-\hat{w})(\infty)=0\text{ order}\geq\theta &}
\end{array}
\right\} .
\end{equation}

\begin{lem} \label{lem:F2:cont}
With the assumptions as above, ${\mathcal
F}=M_{\underline{\theta}}$.
\end{lem}
\begin{Proof}
If $e\in{\mathcal F}$, then $e\in X_{\underline{\theta}}$ as
argued earlier. The constraints on $\hat{e}$ at the $z_i$, $p_i$,
$v_i$ follow as in the discrete--time case. Conversely, if $e\in
M_{\underline{\theta}}$, form $\hat{q}:=
\hat{v}/\hat{n}-\hat{e}/(\hat{n}\hat{w}\hat{d})=\frac{1}{\hat{d}}\left(
\frac{\hat{w}-\hat{e}}{\hat{w}\hat{n}}-\hat{x}\right)$. The only
possible $\overline{\C_+}_e$--poles of $\hat{q}$ are at the
$\overline{\C_+}$--poles/zeros of $P$, the
$\overline{\C_+}$--zeros of $\hat{w}$, and at infinity. The
constraints at each of the $\overline{\C_+}$--points ensure that
$\hat{e}$ has no poles there. It remains to check the behavior at
$\infty$. Now, $\hat{d}(\infty)\ne 0$, and the prescribed behavior
of $\hat{e}$ and $\hat{e}-\hat{w}$ means that
$\frac{\hat{w}-\hat{e}}{\hat{w}\hat{n}}=O(\frac{s^{-\theta_p-\theta_w}}{s^{-\theta_w}s^{-\theta_p}})=O(1)$
for $|s|\to\infty$. Thus $\hat{q}(\infty)$ is finite, and since
$\hat{q}$ has no poles in $\overline{\C_+}$, it is in
$\mathbf{{\mathcal S}}$. It now follows that $e\in{\mathcal F}$.
\end{Proof}

As in the discrete--time case, whenever $z_i$ and $z_j$ form a
conjugate pair, we retain only one of these in the list of
constraints characterizing $M_{\theta}$. We make a similar
reduction for the $p_i$ and $v_i$ also. This entails no loss of
information from $M_{\theta}$ (since $\hat{e}(\bar{z})=
\overline{\hat{e}(z)}$ for any $z$, and any real--valued function
$e$.)

\smallskip

For each interpolation point $z_j=x_j+iy_j$ (recall $x_j>0$), define
$$ \mathbf{a}_j(t):= e^{-x_j t}\cos y_jt\,,\,\,  \mathbf{a}_{j+1}(t):=
e^{-x_j t}\sin  y_jt \,.$$
If $ \mathbf{b}_j$ and $ \mathbf{c}_j$ are also  defined similarly with
 respect to the
$p_j$, $v_j$ respectively,
$M_{\theta}$ takes the form
\begin{equation}
M_{\underline{\theta}}=\left\{e\in X_{\underline{\theta}}\,
\Bigg{|}
\begin{tabular}{l}
$\langle e,\mathbf{a}_i\rangle =\langle w,\mathbf{a}_i\rangle
\quad i=1,2,\ldots ,n$ \\
$ \langle e,\mathbf{b}_i\rangle =\,\,0\qquad\quad \! i=1,2,\ldots ,m$ \\
$ \langle e,\mathbf{c}_i\rangle =\,\,0\qquad\quad i=1,2,\ldots ,l$%
\end{tabular}
\right\} ,  \label{setMtheta}
\end{equation}
where for functions $u\in L^{\infty}$ and $v\in L^1$, $\la u, v \ra:=\int_0^\infty
v(t)u(t)dt$.

Let $A:C_0(\R_+)\to \R^{m+n+l}$ be defined by
 \begin{equation}\label{eq:A}
Ae:= (\langle e,\mathbf{a}_1\rangle,\ldots,\langle e,\mathbf{a}_n\rangle,
\langle e,\mathbf{b}_1\rangle,\ldots,\langle
e,\mathbf{b}_m\rangle, \langle
e,\mathbf{c}_1\rangle,\ldots,\langle e,\mathbf{c}_l\rangle)^T\in
\R^{m+n+l}\,,
\end{equation}
and let
 \begin{equation}\label{eq:b}
b:=(\langle w,\mathbf{a}_1\rangle,
\ldots,\langle w, \mathbf{a}_n\rangle,
0,0,\ldots)^T\in\R^{m+n+l}\,.
\end{equation}
Then, $M_{\underline{\theta}}$ has the form $\{e\in X_{\underline{\theta}}\mid Ae=b\}$.

Given a performance measure
$f:X_{\underline{\theta}}\to \R\cup\{+\infty   \}$,  our question is to evaluate
$$\text{{\bf (P)}:  }\quad\inf_{e\in M_{\underline{\theta}}} f(e)\,,$$
which represents a
theoretical limit of performance for "physically realizable"
controllers (in the sense of  having rational Laplace transform).
This is the central object of study in this paper.
\smallskip

\subsection{primal time-domain performance objectives}\label{sectsometdperf
:cont}

We shall principally consider functionals $f_0$ on $C_0(\R_+)$ of the form:
\begin{eqnarray*}
f_{ma}(e) \!\!\!&:&\!\!\!=\sup_t |e(t) |=\| e\| _\infty \quad\quad\text{(maximum amplitude)}\\
f_{pos}(e) \!\!\!&:&\!\!\!=\sup_t\left[ (e(t))_{+}\right] \qquad\qquad\text{(positive error)}\\
f_{os}(e) \!\!\!&:&\!\!\!=\sup_t\left[ (-e(t))_{+}\right] \quad\qquad\,\,\text{(overshoot)}\\
f_{fl}(e) \!\!\!&:&\!\!\!=\frac12\left[ \sup_t(e(t))-\inf_t(e(t))\right]  \quad\text{(fluctuation)}\\
f_{us}(e)\!\!\! &:&\!\!\!=\sup_t\left[ (e(t)-w(t))_{+}\right] \qquad\quad\,\,\text{(undershoot)}.
\end{eqnarray*}
where for real $\lambda$, we define $\lambda_+:=\max(\lambda,0)$ and
$\lambda_-:=\min(\lambda,0)$. It is straightforward to verify that these functionals are all convex, and also
continuous (in fact Lipschitz, with constant 1) with respect to the $\|\cdot\|_\infty$--norm on $C_0(\R_+)$.

 As in discrete--time, we have $f_{os}+f_{pos}=2f_{fl}$ and
that
\begin{equation}
\left.
\begin{array}{c}
f_{pos} \\
f_{os}
\end{array}
\right\} \leq f_{ma}\leq 2f_{fl}\leq 2f_{ma}.
\label{ineq1:cont}
\end{equation}
Also,  for $e\in C_0(\R_+)$,
\begin{equation}
f_{fl}(e)=\min_{\xi \in {\scriptsize \R} }\sup_t |e(t)-\xi |\, .
  \label{argmin:cont}
\end{equation}
whose proof follows by trivial modification of the proof of its
discrete counterpart in \cite{HEHW}. Later, we will consider some simple time-domain
constraints, so require treatment of functionals of form $f_0+\iota_T$ for
appropriate choices of sets $T$ representing these additional conditions.

\begin{remark}
Recall that a {\em rational\/} error--signal $e(\cdot)$ satisfies (the zero--steady--state condition)  $\lim_{t\to\infty}e(t)=0$ iff $e\in C_0(\R_+)$
iff $e\in L^p(\R_+)$ ($p\ne\infty$). Thus, we {\em could\/} have chosen   $L^p$ ($p\ne\infty$)  as our ambient space---%
however, the functionals $f_0$ above are not continuous relative to these $L^p$--norms, negating the usefulness of this choice.
 \end{remark}

\subsection{initial statement of duality for our primal problem}\label{sec:dual:statement}
The analysis of problem {\bf (P)} will proceed by recasting it in dual form via Proposition~\ref{prop:BL}, as  is given
below in Proposition~\ref{thmfench:cont}. As preparation, we require the following characterisation of the
closure $\overline{X_{\underline{\theta}}}$ (proved in the Appendix).
 \begin{prop}\label{prop:summary}\,\,\,\quad $\overline{X_{(\theta_w,\theta_p)}}=
\left\{
\begin{array}{ll}C_{0,\alpha}(\R_+)\,\,(\text{with }\alpha=w(0+)) & \,\mbox{ if }\,\,\theta_p>0,\,\theta_w=1\\
C_{0}(\R_+) & \,\mbox{ if }\,\,\theta_p=0,\,\theta_w=1\\
C_{0,0}(\R_+) & \,\mbox{ if }\,\,\theta_w>1
\end{array}
\right.$
(where $C_{0,\alpha}(\R_+)\subseteq C_0(\R_+)$ denotes the  set of functions $\varphi$ for which $\varphi(0)=\alpha$).
\end{prop}

The formulation of the dual problem will also require the following spaces:
$$U:=\func{span}[ \mathbf{a}_1,\ldots , \mathbf{a}_n],\
V:=\func{span}[ \mathbf{b}_1,\ldots , \mathbf{b}_m],\
W:=\func{span}[ \mathbf{c}_1,\ldots , \mathbf{c}_l],\,$$
 From the blanket assumptions on poles and zeros, we have
$U,\ V,\ W$ contained in $C_0(\R_+)$. Also, from the resulting integrability of the $\mathbf{a}_n$, $\mathbf{b}_n$, $\mathbf{c}_n$ , each
 $e^*\in U\oplus
V\oplus W$ defines a measure $\mu\in \mathbf{M}(\R_+)$ by
\begin{equation}\label{eq:meas}
 \mu (E)=\int_E e^* \qquad\text{ for Borel sets $E$}\,.
\end{equation}
Thus $U\oplus V\oplus W $ may also be considered as a subspace of
$X^*={\mathbf M}(\R_+)$.

The next basic duality result is the foundation for all subsequent analysis.
 \begin{prop}
\label{thmfench:cont}
Let $T\subseteq C_0(\R_+)$ be convex, $f_0:C_0(\R_+)\to \R$ convex and continuous, with
$b\in A(X_{\underline{\theta}}\cap\func{int} T)$, where $A$ and $b$ are given in (\ref{eq:A}) and
(\ref{eq:b}).

If $\theta_w + \theta_p=1$ (i.e.\ $\overline{X_{\underline{\theta}}}=C_0(\R_+)$ via Prop~\ref{prop:summary}) then
\begin{equation}\label{eq:duality:1}
\inf_{ T\cap X_{\underline{\theta}}\cap A^{-1}b }\ f_0  =\max_{
{\scriptsize
\begin{array}{c}
 \mu\in U\oplus V\oplus W \\
%\mu\in\func{dom}f^{*} \\
\end{array}
} } \left[ \langle \func{Proj}_U(\mu),w\rangle -(f_0 +\iota_T)^{*}(\mu)\right]
,
\end{equation}
and if instead, $\theta_w + \theta_p>1$ (so that by Prop~\ref{prop:summary},  $\overline{X_{\underline{\theta}}}=C_{0,\alpha}(\R_+)$ with $\alpha=w(0+)$) then
\begin{equation}\label{eq:duality:2}
\inf_{ T\cap X_{\underline{\theta}}\cap A^{-1}b }\ f_0  =\max_{
{\scriptsize
\begin{array}{c}
 \mu\in U\oplus V\oplus W \\
%\mu\in\func{dom}f^{*} \\
\end{array}
} } \left[ \langle \func{Proj}_U(\mu),w\rangle +(f_0 +\iota_T)^{\#}(\mu)\right]
,
\end{equation}
where: $\func{Proj}_U (\cdot)$ denotes the natural projector  from
$U\oplus V\oplus W$ onto $U$; and, for any $\mu$, and any $f$,
\begin{equation}\label{eq:duality:3}
f^{\#}(\mu):= \max_ {\lambda\in{\scriptsize
\R}}\,[\alpha\lambda-f^*(\mu+\lambda\delta)]
\end{equation}
with $\delta$ denoting the Dirac measure at $0\in\R_+$.

%
%
% If the infimum on the left is achieved for some $e_{opt}\in
%M$, then
%\begin{equation} \label{align:cont}
%f(e_{opt})+f^{*}(\mu_{opt})=\langle e_{opt},\mu_{opt}\rangle .
%\end{equation}
\end{prop}
\begin{Proof}
Now, $b\in A(X_{\underline{\theta}}\cap\func{int} T)$ implies, by convexity, that
$\overline{T\cap X_{\underline{\theta}}\cap A^{-1}b}\supseteq T\cap \overline{  X_{\underline{\theta}}\cap A^{-1}b}$.
But $\overline{  X_{\underline{\theta}}\cap A^{-1}b}=\overline{  X_{\underline{\theta}}}\cap A^{-1}b $ by a simple argument using
the finite-codimensionality of $A^{-1}b$. Thus $\overline{T\cap X_{\underline{\theta}}\cap A^{-1}b}\supseteq
T\cap \overline{  X_{\underline{\theta}}}\cap A^{-1}b$ and hence by continuity of $f_0$,
$$
\inf_{ T\cap X_{\underline{\theta}}\cap A^{-1}b }\ f_0 =\inf_{ T\cap \overline{X_{\underline{\theta}}}\cap A^{-1}b }\ f_0 \,.
$$
We now apply Proposition~\ref{prop:BL} to the latter minimization in $\overline{X_{\underline{\theta}}}$. Place $f:=f_0+\iota_T$
 and $g:=\iota_{\{b\}}$.
If $\overline{X_{\underline{\theta}}}=C_0$, and taking  $X=C_0$, then $b\in A(\func{int dom}f)$ and after some simple manipulations of
the resulting dual, we obtain the form
 (\ref{eq:duality:1}).

 If instead, $\overline{X_{\underline{\theta}}}=C_{0,\alpha}$, note that
 $ A^{-1}b \cap\overline{  X_{\underline{\theta}}}=\widetilde{A}^{-1}\widetilde{b}$, where
 $\widetilde{A}e:=(Ae,e(0))^T$ (for all $e\in C_0(\R_+)$) and $\widetilde{b}:=(b,\alpha)^T$.
 We also have $\widetilde{b}\in \widetilde{A}(\func{int}T)$, so we may apply the duality to
 $\inf_{ T\cap C_{0,\alpha}\cap A^{-1}b }\ f_0 =\inf_{ T\cap   \widetilde{A}^{-1}\widetilde{b} }\ f_0 $,
 with $X=C_0$ again, with the new $\widetilde{A}$, $\widetilde{b}$, which eventually yields the dual in
 (\ref{eq:duality:2}).
\end{Proof}

\section{Dual Formulation}
\label{dualformulation:cont}
The dual characterisation will be completed by evaluation of the conjugate functionals
appearing on the right-hand-side of equations (\ref{eq:duality:1}) and (\ref{eq:duality:2}) above.
We now study the forms of the conjugates for a general class of
cost--functions which includes those of interest in this paper.
%A simple representation for the domains of the conjugates is
%derived, from which the dual constraints can be written down
%directly.

We assume the objective functionals $f:C_0(\R_+)\to \R\cup\{+\infty\}$ take the form
\begin{equation}\label{eq:fF}
 f(e)=\sup_{t\ge 0}F_t(e(t))
\end{equation}
where $F$ satisfies the following assumptions.
\begin{enumerate}
\item\label{ass:F:1}
  For each $t\geq 0,\, F_t:\R\to [0,+\infty]$
%\,\,\,\, 0\in
%F(\R);\,\,\,\,\exists L>0\text{ with }0\in [e^{-}(L),e^{+}(L)]\,,
 \item\label{ass:F:2}
For all $t\ge 0$, all $L\ge 0$, the sublevel--set $[F_t\le L]$ is a nonempty closed (possibly unbounded) interval in $\R$, with endpoints
$e_t^{-}(L)\in \R\cup\{ -\infty\}$, $e_t^{+}(L)\in  \R\cup\{ +\infty\}$ respectively, with, further,
\[
\inf_{t\ge 0}e_t^{+}(L)>-\infty\,,\quad \sup_{t\ge 0}e_t^{-}(L)<+\infty\,.
\]
 \item\label{ass:F:3}
 For each $L\ge 0$, the functions $t\mapsto e_t^\pm (L)$ are piecewise continuous
 (in appropriate sense for extended-real-valued functions)
 with at worst a countable set of jump-discontinuities, where at all such jumps $\bar{t}$, have
 \[
 \min\{ e_{\bar{t}+}^+(L), e_{\bar{t}-}^+(L)    \} \ge  \max \{e_{\bar{t}+}^-(L), e_{\bar{t}-}^-(L)     \}\,.
 \]
 \item \label{ass:F:4} For each $L\ge 0$, there is $t_0(L)\ge 0$ such that for all $t\ge t_0$, have that
 $e_t^-(L)\le 0\le e_t^+(L)$.
\end{enumerate}

 Before we proceed further, a review of some relevant measure
theory is appropriate. (Our reference shall be \cite{Cohn}.)
 Given $\mu\in\mathbf{M}(\R_+)$, we can find a
{\it Hahn decomposition} of $\R_+$ into disjoint Borel sets $P$,
$N$ such that $\mu(E\cap P)\ge 0$ and $\mu(E\cap N)\le 0$ for each
Borel set $E$. From this are obtained regular finite Borel
measures $\mu_{\pm}$ given by $\mu_+(E):=\mu(E\cap P)$ and
$\mu_-(E):=\mu(E\cap N)$, with $\mu_+\ge 0$ and $\mu_-\le 0$. We
then have $\mu=\mu_++\mu_-$ and the {\it variation} of $\mu$
(denoted $|\mu|$) is defined by $|\mu|:=\mu_+-\mu_-$. Note that
our sign convention has been chosen to conform with that used in
the discrete-time analysis of \cite{HEHW}, but differs from the
usual choice in measure theory, where $\mu_+$ and $\mu_-$ are both
non--negative, whereas here we have $\mu_-\le 0$. Recall that we
assume $\mathbf{M}(\R_+)$ to be normed by $\|\mu\|:=|\mu|(\R_+)$.

Note that if measure $\mu\in U\oplus V\oplus W$ and function $e^*$
are related via \eqref{eq:meas}, then by standard arguments follows that
$$
\mu_\pm (E)=\int_E e_\pm ^* \,\quad\text{ for Borel sets $E$}\,,
$$
 and that $\mu\ge 0$ if and only if $e^*(t)\ge 0$ for all $t$, with a similar
relation for the reverse inequality $\mu\le 0$.

We also recall the following standard measure-theoretic convention: Since we need to integrate
real functions taking $+\infty$ as a possible value, the definition of the integral
(with respect to a measure) incorporates the convention
$0\cdot\infty:=0$ whenever one of the factors is the value of a (unsigned) measure.

\subsection{The Conjugate for the Assumed Form of the Primal Objective}
The  following result provides a  continuous--time analogue
of part of \cite[Theorem 9]{HEHW}. (The proof is deferred to the Appendix.)

\begin{thm}\label{thmfstar1:cont}
For $f:C_0(\R_+)\to \R\cup\{+\infty\}$  defined as in
\eqref{eq:fF}, where $F$ satisfies the assumptions 1---4 above. Then,

\begin{equation}\label{eq:conj:cont}
f^*(\mu)=\sup_{  L\ge 0  } \left(\int_0^{+\infty}\!\! \!e_t^+(L)d\mu_+(t)+
 \int_0^{+\infty }\!\!\!e_t^-(L)d\mu_-(t)-L\right)
\end{equation}
for all $\mu\in \mathbf{M}(\R_+)$ such that $|\mu|(\{\bar{t}\})=0$ for each $\bar{t}$ for which there is some $L$
for which at least one of $e^+(L)$ or $e^+(L)$ is discontinuous at $\bar{t}$.
 \end{thm}

\begin{remark}\label{rem:thmfstar1:cont}
 If $F_t\equiv F$ (independent of $t$), assumption 3.\ is inactive, and assumption 4.\ can only be valid for $L$ satisfying
 $L\ge F(0)$ (which is equivalent to $[f\le L]\neq\emptyset$) and from the proof (see Appendix) it follows that the supremum in
 (\ref{eq:conj:cont}) is  then to be taken over $L\ge F(0)$.
\end{remark}

\subsection{Duals for some Time-Domain Minimization Problems}
We may now derive the duals of the continuous--time versions of
MA, OS, POS, FL in the manner of \cite{HEHW}.

Because of Proposition~\ref{prop:summary}, this separates into
cases where $\overline{X_\theta}=C_0(\R_+)$ and
$\overline{X_\theta}=C_{0,\alpha}(\R_+)$. We begin with the former case.

\subsubsection{Maximum Amplitude}
Clearly
$e^{+}(L)=L$ and $e^{-}(L)=-L.$ By Theorem~\ref{thmfstar1:cont} and Remark~\ref{rem:thmfstar1:cont},
$$\func{dom}%
f_{ma}^{*}=\left\{ \mu\in \mathbf{M}(\R_+)\mid\| \mu\|\leq 1\right\}
\text{ and }f_{ma}^{*}(\mu)=0\,.$$

Thus by Proposition~\ref{thmfench:cont}, the dual of the problem of
maximum--amplitude
minimization, denoted MA DUAL, can be written as
\begin{equation}
\max_{{\scriptsize \begin{array}{c}  \mu\in U\oplus V\oplus W  \\
 \|\mu \|\leq 1
\end{array}}}  \langle \func{Proj}_U(\mu),w\rangle .
\label{madual:cont}
\end{equation}

\subsubsection{Positive Error}

For this case we have $e^{-}(L)=-\infty$   and $e^{+}(L)=L$. By Theorem~\ref%
{thmfstar1:cont} and Remark~\ref{rem:thmfstar1:cont}, $\func{dom}f_{pos}^{*}=\left\{ \mu\in
\mathbf{M}(\R_+)\mid \mu\geq 0%
\text{ and }\mu(\R_+)\leq 1\right\} $ and $%
f_{pos}^{*}(\mu)=0$. The dual of POS is
\begin{equation}
\max_{{\scriptsize \begin{array}{c}  \mu\in U\oplus V\oplus W
  \\ \mu\geq 0\text{, }\|%
\mu\|\leq 1  \end{array}}}  \langle \func{Proj}_U(\mu),w\rangle ,
\label{posdual:cont}
\end{equation}

\subsubsection{%
Overshoot} For overshoot minimization $e^{+}(L)=+\infty$ and
$e^{-}(L)=-L$. By Theorem~\ref{thmfstar1:cont} (and the Remark),
$\func{dom}f_{os}^{*}=\left\{ \mu\in \mathbf{M}(\R_+)\mid \mu\leq
0\text{ and }\mu(\R_+)\geq -1\right\} $
and $%
f_{os}^{*}(\mu)=0.$ The dual of OS is
\begin{equation}
\max_{{\scriptsize \begin{array}{c}  \mu\in U\oplus V\oplus W
  \\ \mu\leq 0\text{, }\|%
\mu\|\leq 1  \end{array}}}  \langle \func{Proj}_U(\mu),w\rangle
\text{.}  \label{osdual:cont}
\end{equation}

\subsubsection{Fluctuation}
Also, as in the discrete--time case, we may use
\eqref{argmin:cont} to deduce that $f_{fl}^*$ is the indicator
function of
$$ \{\mu\in\mathbf{M}(\R_+)\mid \mu_+(\R_+)\le \frac12,\,
 \mu_-(\R_+)\ge -\frac12\}\,,$$
(see proof of \cite[Theorem 11]{HEHW}). Hence, by
Proposition~\ref{thmfench:cont}, the dual of FL is
$$ \max_{\scriptsize \begin{array}{c} \mu\in U\oplus V\oplus W\\
                             \mu_+(\R_+)\le\frac12\\
                           \mu_-(\R_+)\ge -\frac12 \end{array} }
           \la \func{Proj}_U (\mu),w\rangle   \,. $$

We observe in passing that for the fluctuation--minimization problem,
that the minimum is never achieved (even over $\overline{M_{\underline{\theta}}}$),
except in the trivial case where $P$ has no
poles in $\overline{\C_+}$, so $0\in M:=\overline{M_{\underline{\theta}}}$ is optimal.
Indeed, suppose an optimal $e\in M$ is attained, and let $\mu$ denote the optimal dual
element. These must satisfy the
 alignment condition $ f_{fl}(e)+f_{fl}^*(\mu)= \la e,\mu\ra  $,
  which yields $$ \la e,\mu\ra= f_{fl}(e)=\frac12 (\sup e-\inf e)$$
since $f_{fl}^*(\mu)=0$. If $\text{FL}_{opt}=0$ this immediately
yields the contradiction $e=0$ (since $0\notin M$). If
$\text{FL}_{opt}$ is positive, note that the optimal $\mu$ is then
nonzero.
 Now as $\int d\mu_+\le\frac12$ and
$\int d\mu_-\ge -\frac12$, and $\sup e\ge 0\ge\inf e$, we obtain
\begin{align*}
\int e\,d\mu_+ + \int e\,d\mu_-=\la e,\mu\ra   &= \frac12 (\sup e-\inf e)\\
& \ge \sup e \int d\mu_+ + \inf e \int d\mu_- \,,
\end{align*}
so that
$$ 0\le \int (e-\sup e)\,d\mu_+ +\int (e-\inf e)\,d\mu_- \le 0\,.$$
Since both terms in the above are nonpositive, we conclude that
$$ \int (e-\sup e)\,d\mu_+ = \int (e-\inf e)\,d\mu_-=0\,.$$
Writing this in terms of the function $e^*$ associated to $\mu$ via
\eqref{eq:meas}, we have for almost all $t$ (w.r.t.\ Lebesgue measure),
and hence by continuity, for all $t$, that
$$ (e(t)-\sup e)e^*(t)=0=(e(t)-\inf e)e^*(t)\,. $$
Thus
$$ e(t)=\left\{ \begin{array}{cl}
 \sup e & \text{if }e^*(t)>0\\
 \inf e &\text{if }e^*(t)<0
\end{array} \right. .$$
By analyticity of $e^*\ne 0$ (being a finite sum of sinusoids), each of
its zeros is isolated, so by continuity of $e$ follows that $\inf
e=\sup e$, implying  again the contradiction $e=0$.

\smallskip
\subsubsection{Undershoot}
 We assume $w(t)\ge 0$ for all $t$. (Note: To derive an expression for
 $f^*_{us}$,
  we do not yet need to assume a rational
Laplace transform for $w$). Place
$F_t(\xi):=(\xi-w(t))_+$. Then $F_t(0)=0$ for all $t\ge 0$, and
$0\in\{F_t\le L\}=(-\infty, w(t)+L]$ for each $L\ge 0$.   Theorem~\ref{thmfstar1:cont}
gives,  for all $\mu\in {\bf M}(\R_+)$, that
$$
f^*_{us}(\mu)=\int w\,d\mu_+ + (-\infty)\cdot\mu_-(\R_+) +\sup_{L\ge 0} L(\mu_+(\R_+)-1)
$$
which yields the following
\begin{cor}\label{prop:undershoot:conj}
Let $w\ge 0$.  Then,
$$
\func{dom} f^*_{us}\subseteq \{ \mu\in {\bf M}(\R_+)\mid \mu\ge
0\}
$$
and for each $\mu\ge 0$,
$$
f^*_{us}(\mu)=\left\{ \begin{array}{cl}
 \int w \,d\mu & \text{if }\mu\ge 0\text{ and } \|\mu\|\le 1\\
 +\infty &\text{otherwise }
\end{array} \right.\,.
$$
Note that $\int w \,d\mu $ may take an infinite value.
\end{cor}

Consequently, the dual for US has the form
\begin{equation}
\max_{{\scriptsize \begin{array}{c}  \mu\in U\oplus V\oplus W
  \\ \mu\leq 0\text{, }\|%
\mu\|\leq 1  \end{array}}}  \langle \func{Proj}_{V\oplus
W}(\mu),w\rangle \text{.}  \label{osdual:cont}
\end{equation}

\smallskip

\bigskip

\smallskip

Recall that each $\mu\in U\oplus V\oplus W$ corresponds to a
function $e^*$ by Equation~\eqref{eq:meas}. This has the
consequences: $\|\mu\|=|\mu|(\R_+)=\int |e^*|=\|e^*\|_1$;  $\mu\ge
0$ (as a measure) if and only if $e^*(t)\ge 0$ for all $t$; and
$\mu\le 0$ if and only if $e^*(t)\le 0$ for all $t$. We may
therefore restate these duals in  a form identical to those in
discrete time:

\begin{align*}
 \text{(MA DUAL)}\qquad &
\max_{{\scriptsize \begin{array}{c}  e^*\in U\oplus V\oplus W  \\
 \|e^* \|_1\leq 1
\end{array}}}  \langle \func{Proj}_U(e^*),w\rangle \\
 \text{(POS DUAL)}\qquad  &
\max_{{\scriptsize \begin{array}{c}  e^*\in U\oplus V\oplus W
  \\ e^*\geq 0\text{, }\|%
e^*\|_1\leq 1  \end{array}}}  \langle \func{Proj}_U(e^*),w\rangle   \\
\text{(OS DUAL)}\qquad &
\max_{{\scriptsize \begin{array}{c}  e^*\in U\oplus V\oplus W
  \\ e^*\leq 0\text{, }\|%
e^*\|_1\leq 1  \end{array}}}  \langle \func{Proj}_U(e^*),w\rangle \text{} \\
\text{(US DUAL)}\qquad & \max_{{\scriptsize \begin{array}{c}
e^*\in U\oplus V\oplus W
  \\ e^*\leq 0\text{, }\|%
e^*\|_1\leq 1  \end{array}}}  \langle \func{Proj}_{V\oplus W}(e^*),w\rangle \text{} \\
\text{(FL DUAL)}\qquad &
 \max_{\scriptsize \begin{array}{c} e^*\in U\oplus V\oplus W\\
                            \int_{\tiny  \R_+}e_+^*\le\frac12\\
                            \int_{\tiny  \R_+}e_-^*\ge -\frac12
\end{array} }
           \la \func{Proj}_U (e^*),w\rangle  \,.
\end{align*}

\subsubsection{case of  $\overline{X_\theta}=C_{0,\alpha}$   }
When
$\overline{X_\theta}=C_{0,\alpha}$, we apply the duality formula
(\ref{eq:duality:2}), which requires calculation of $f^{\#}$,
given by (\ref{eq:duality:3}). The forms of all the duals, except
for FL, will be unchanged, because of the next result, proved in
the Appendix.
\begin{lem}\label{lem:fhash}
If $f$ stands for any of the functionals treated (except for
$f_{fl}$), then for $\mu\in U\oplus V\oplus W$, have
$$f^{\#}(\mu)=-f^*(\mu) \,, $$
whereas for FL, we have
$$
f_{fl}^{\#}(\mu)=-f_{fl}^*(\mu)+\alpha(
1/2-\mu_+(\R_+))\,.\qquad\qquad \square
$$
\end{lem}
%\newline\noindent
This, for FL, yields the dual
$$
\text{(FL DUAL)}\qquad
 \max_{\scriptsize \begin{array}{c} e^*\in U\oplus V\oplus W\\
                            \int_{\tiny  \R_+}e_+^*\le\frac12\\
                            \int_{\tiny  \R_+}e_-^*\ge -\frac12
\end{array} }
           \la \func{Proj}_U (e^*),w\rangle +\alpha(1/2-\int_{\tiny  \R_+} e_+^*) \,.
$$

\subsection{Overshoot/Undershoot Minimization}\label{sectosusmin:cont}

We follow  the procedure of \cite{HEHW}. The first lemma is
essentially the continuous--time counterpart of \cite[Proposition
15]{HEHW}.

\begin{lem}\label{prop:sinus:cont}
Let $y_1,\ldots, y_N\ne 0$ with $y_i\ne\pm y_j$ whenever $i\ne j$. If
$\alpha_1,\ldots,\alpha_N$  and $\beta_1,\ldots,\beta_N$ are real scalars
such that for some $C,\ \rho>0$ we have
$$ \sum_{i=1}^N (\alpha_i \cos y_it +\beta_i \sin y_it)\ge -Ce^{-\rho t}
\text{ for all $t$ large} $$
then $\alpha_i=0=\beta_i$ for $i=1,\ldots,N$.
\end{lem}

\begin{Proof}
Place $a(t):= \sum_{i=1}^N (\alpha_i \cos y_it +\beta_i \sin y_it)$.
For $r>0$, set $$S(r):=\int_0^\infty a(t)e^{-rt}dt=
\sum_{i=1}^N \frac{\alpha_ir+\beta_i y_i}{r^2+y_i^2}\,.$$
Since $y_i\ne 0$ for all $i$, there exists $\lambda>0$ such that
$|S(r)|\le \lambda$ for all $r>0$.
Thus for  some $t_0$ and any $N\in\R$ and $r>0$,
\begin{align*}
0\le \int_0^N a_+(t)e^{-rt}dt &\le \int_0^\infty  a_+(t)e^{-rt}dt \\
& = S(r)-\int_0^\infty  a_-(t)e^{-rt}dt\\
& \le \lambda + \int_0^{t_0}|a_-(t)|dt+\int_{t_0}^\infty |a_-(t)|dt \\
&\le \lambda +\int_0^{t_0}|a_-(t)|dt+ C\int_0^\infty e^{-\rho t}dt\\
&\le \lambda +\int_0^{t_0}|a_-(t)|dt+ C/\rho :=\lambda '
\end{align*}
For each fixed $N$, we may let $r\to 0$ in the above, to obtain
by the Dominated Convergence Theorem that
$0\le\int_0^N a_+(t)dt \le \lambda '$. Since $|a_-(t)|\le Ce^{-\rho t}$ for all
large $t$, it follows that $\int_0^N|a(t)|dt$ is bounded above for all $N$,
so that $a(\cdot)$ is in $L^1$. This implies that the Laplace transform
$\hat{a}(s)$ is defined for all $s$ such that $\func{Re}s\ge 0$, and
is bounded in this region. However, for all $s$ with positive real part,
$$ \hat{a}(s)= \sum_{i=1}^N\frac{\alpha_is+\beta_i y_i}{s^2+y_i^2} $$
and so if $(\alpha_k, \beta_k)\ne (0,0)$, then from the assumed distinctness
of the $y_i$, follows that $\hat{a}$ must have a pole at $s=\pm iy_k$,
contradicting the boundedness of the transform at points of the imaginary
axis. Thus necessarily $(\alpha_k, \beta_k)=(0,0)$.
\end{Proof}

\smallskip

Lemma~\ref{prop:sinus:cont} can now yield an analogue of
\cite[Proposition 17]{HEHW} in continuous time. Given the complex
numbers $z_1,.., p_1,..,$ and $v_1,..$, let $\gamma$ denote the
smallest value among those poles/zeros that lie on the positive
real axis; if there are none, set $\gamma$ to $\infty$.
\begin{thm}\label{cor1:cont}
Consider the `equivalence' class $\mathcal{A}$, of all plants
which: \begin{enumerate}
\item  have the same $\gamma$ value; and
\item  have the same poles and zeros in the region $\{s\mid \func{Re}s\ge
\gamma\}$.
\end{enumerate}
Then, for any fixed minimum--phase reference input
(i.e.\ having no zeros in the closed right half--plane), all plants
 in $\mathcal{A}
$ have the same value for $\funcc{OS}_{\scriptsize{\funcc{opt}}}$. If the
reference--input is also non-negative, then all plants
 in $\mathcal{A}
$ have the same value for $\funcc{US}_{\scriptsize{\funcc{opt}}}$.
\end{thm}
\begin{Proof}
We show that the poles/zeros in $\{s\mid\func{Re}s<\gamma\}$ do
not affect the value of the dual.

For simplicity, we assume that all the
  poles/zeros in $\{s\mid\func{Re}s<\gamma\}$ have equal
real part. The general case is handled by iteration of the argument given below.
Let $e^*\in U\oplus V\oplus W$ with $e^*\ge 0$ (or $e^*\le 0$). Regrouping
the modes of $e^*$ according to the size of the real parts of the associated
poles and zeros gives
$$ e^*(t)=e^{-\rho_1 t}
\sum_{i=1}^{N_1}(\alpha _i^{(1)}\sin (y_i^{(1)}t)+
\beta _i^{(1)}\cos (y_i^{(1)}t)) +
e^{-\rho_2 t}\sum_{i=1}^{N_2}(\alpha _i^{(2)}\sin (y_i^{(2)}t)+
\beta _i^{(2)}\cos (y_i^{(2)}t)) + \ldots
$$
where $0<\rho_1<\rho_2<...$, and in the first sum, which corresponds to the
 poles/zeros in $\{s\mid\func{Re}s<\gamma\}$, we have
$y_i^{(1)}\ne 0$, and $y_i^{(1)}\ne\pm y_j^{(1)}$ whenever $i\ne j$.
Hence
\begin{align*}
\sum_{i=1}^{N_1}(\alpha _i^{(1)}\sin (y_i^{(1)}t)+
\beta _i^{(1)}\cos (y_i^{(1)}t)) &= e^*(t)/e^{-\rho_1 t} +
e^{(\rho_1-\rho_2)t}\sum_{i=1}^{N_2}\ldots \\
& \ge 0 + O(e^{(\rho_1-\rho_2)t})\,.
\end{align*}
By Lemma~\ref{prop:sinus:cont}, $\alpha _i^{(1)}\!=0=\beta
_i^{(1)}$ for $i=1,..,N_1$. Thus the modes associated with  poles
or zeros with real part smaller than $\gamma$ do not contribute to
the dual, as claimed.
\end{Proof}

%TCIMACRO{\FRAME{ftbpF}{8.7778in}{4.683in}{0pt}{}{}{osplot2.eps}%
%width 8.7778in;  height 4.683in;  depth 0pt;  original-width 13.3389in;
%original-height 6.0554in;  cropleft "0.0978";  croptop "1";
%cropright "0.9510";  cropbottom "0";
%filename 'osplot2.eps';file-properties "XNPEU";}}}%

\bigskip

\bigskip

\subsection{Analytical Results for a First--Order Plant}

We first derive analytical expressions
for performance limitations in terms of pole/zero locations for first-order plants.
 Thanks to Theorem~\ref{cor1:cont},
 in the case of overshoot and undershoot limitations on performance, these analytical results
 can be extended to include cases where the plant has an arbitrary number of oscillatory poles
 or zeros (i.e.\ those off the real axis), and the reference input has an arbitrary number of oscillatory zeros.

%\subsubsection{Overshoot, Maximum Amplitude, and Fluctuation}

\begin{prop}
\label{propcontex}Suppose the plant, $P=\frac{s-z_{1}}{s-p_{1}},$
has one real positive zero, $z_{1,}$ and one real positive pole,
$p_{1},$ where $z_{1}>p_{1}.$ Let $\widehat{w}(s)=1/s$. Define
$h:=p_{1}/(z_{1}-p_{1}).$ Then
\begin{align*}
(i)\,\,\qquad \funcc{OS}_{\scriptsize{\funcc{opt}}}  & \, =\,h\\
(ii) \qquad \funcc{MA}_{\scriptsize{\funcc{opt}}}  &  \,=\,\frac1{1-2^{-1/h}}\\
(iii)\,\qquad \funcc{FL}_{\scriptsize{\funcc{opt}}}  & \, =\,\frac{(h+1)^{(h+1)}}{2h^{h}}\\
\!(iv) \,\,\quad \funcc{POS}_{\scriptsize{\funcc{opt}}}  & \, =\, 1.
\end{align*}

\end{prop}

\begin{remark}
It can be verified using elementary calculus that, for $h\in(0,\infty),$%
\[
\max(1,h)\leq\frac1{1-2^{-1/h}}\leq\frac{(h+1)^{(h+1)}}{h^{h}}\leq
\frac2{1-2^{-1/h}}
\]
 From (ii) and (iii) it follows that
$2\funcc{MA}_{\scriptsize{\funcc{opt}}}=\funcc{FL}_{\scriptsize{\funcc{opt}}}$ if and only if $h=1$, that is
$z_{1}=2p_{1}.$ If $h\neq 1$, then $\funcc{FL}_{\scriptsize{\funcc{opt}}}$ is strictly
less than $2\funcc{MA}_{\scriptsize{\funcc{opt}}}$; in this example, when $h\neq 1$, a
minimal fluctuation response is obtainable only by allowing the
infinity--norm to be larger than optimal.
\end{remark}

\begin{Proof}
We have $A=\operatorname*{span}\left[  e^{-z_{1}t}\right] ,\
B=\operatorname*{span}\left[  e^{-p_{1}t}\right]  ,\ C$ is empty
and $w=\underline{1}$ (the unit step function).

(i) The dual optimal vector for (OS DUAL),
$e^{*}=\zeta(t)-\eta(t)$ where $\zeta\in A$ and $\eta\in B,$ will
satisfy $\zeta(0)=\eta\left(  0\right)  ,$ so $\zeta(t)=\alpha
e^{-z_{1}t}$ for some real number $\alpha$ and
$\eta(t)=\alpha e^{-p_{1}t}.$ Hence $\alpha\int_{0}^{\infty}(e^{-p_{1}%
t}-e^{-z_{1}t})dt\leq1,$ from which $\alpha\leq
p_{1}z_{1}/(z_{1}-p_{1}).$ Then $\funcc{OS}_{\scriptsize{\funcc{opt}}}=\max_{\zeta\in
A}\langle\zeta,r\rangle=\max_{\alpha}\int
_{0}^{\infty}\zeta(t)\underline{1}dt=\max_{\alpha}\left[
\alpha/z_{1}\right] =h.$

(ii) (MA DUAL) is
\begin{align*}
&  \max_{\alpha,\beta}\left[  \alpha/z_{1}\right] \\
\text{subject to }  &   \int_0^\infty |\alpha e^{-z_{1}t}-\beta e^{-p_{1}t}|\,dt%
 \leq1.
\end{align*}

It is clear that, at optimality, $\alpha>\beta>0,$ and there will
exist a positive number $t_{0}$ such that
$\zeta(t_{0})=\eta(t_{0}).$ It is obvious also that at optimality
the constraint inequality will be satisfied as an equality. The
dual becomes
\begin{align*}
&  \max_{\alpha,\beta,t_{0}}\left[  \alpha/z_{1}\right] \\
\text{s.t. }\alpha e^{-z_{1}t}  &  =\beta e^{-p_{1}t}\text{ and
}\int _{0}^{t_{0}}\left(  \alpha e^{-z_{1}t}-\beta
e^{-p_{1}t}\right) dt+\int_{t_{0}}^{\infty}\left(  \beta
e^{-p_{1}t}-\alpha e^{-z_{1}t}\right) dt=1.
\end{align*}

This is a simple finite-dimensional constrained optimization
problem, which can be solved using the method of Lagrange
multipliers. After some elementary algebra one obtains
$t_{0}=\frac{\log(\alpha-\beta)}{z_{1}-p_{1}}=-\frac
{\log1/2}{p_{1}},$ $\alpha/\beta=2^{h}$ and $\alpha=z_{1}/[1-2(1/2)^{z_{1}%
/p_{1}}]$, from which the result follows.

(iii) At optimality the inequalities in the two constraint
equations for (FL DUAL) will be satisfied as equalities. Then (FL
DUAL) becomes
\begin{align*}
&  \max_{\alpha,\beta}\left[  \alpha/z_{1}\right] \\
\text{s.t. }  &  \int_0^\infty |\alpha e^{-z_{1}t}-\beta
e^{-p_{1}t}|\,dt =1\text{ and }\int_{0}^{\infty}\left(
\alpha e^{-z_{1}t}-\beta e^{-p_{1}t}\right)  dt=0.
\end{align*}

The second constraint gives immediately that $\alpha p_{1}=\beta
z_{1}.$ Let
$t_{0} $ be the positive number with the property that $\zeta(t_{0}%
)=\eta(t_{0}).$ After performing the integration in the first
constraint, and writing it as an equality, the dual becomes
\begin{align*}
&  \max_{\alpha,t_{0}}\left[  \alpha/z_{1}\right] \\
\text{s.t. }1  &  =\frac{2\alpha}{z_{1}}\left(  e^{-p_{1}t_{0}}-e^{-z_{1}%
t_{0}}\right)  .
\end{align*}

After some algebra the optimizing $t_{0}$ and $\alpha$ are found
to be
$t_{0}=\frac{\log(z_{1}-p_{1})}{z_{1}-p_{1}}$ and $\alpha=z_{1}(1+1/h)^{h}%
(1+h)/2.$

(iv) In (POS DUAL), if $r\geq0,$ an optimizing $\eta$ will be
identically zero. Then the optimizing $\zeta$ will be positive and
satisfy $\parallel \zeta\parallel_{1}=1.$ For $w=$\underline{$1$}
the optimal dual cost is
$\int_{0}^{\infty}\zeta$\underline{$1$}$dt=1.$\quad
\end{Proof}

\bigskip

\bigskip

\subsection{A time-domain constraint effect}

We consider some relationships between overshoot (and undershoot)
with some simple finite-time-horizon constraints. (We have in mind the effect
of rise-time constraints on optimal overshoot performance.)
Time-domain signal bounds will be represented by the set
  \begin{equation}%\label{eq:T}
  T:=\{e\in C_0(\R_+)\mid \phi_-(t)\le e(t)\le \phi_+(t) \,\,\forall t\ge 0\}
 \end{equation}
 for suitable bounding functions $\phi_{\pm}:\R_+\to \R\cup \{ -\infty, +\infty  \}$.
For simplicity, we restrict to the finite-horizon case, where for some positive
$\bar{t}>0$, have $\phi_{\pm}(t)\equiv \pm\infty$ for all $t\ge \bar{t}$.
Also, assume that $\phi_{\pm}$ are finite-valued and continuous on $[0,\bar{t}]$,
with $\phi_{-}(t) < \phi_{+}(t)$ for all $t$, and $\phi_{+}(t)\ge 0$. In this case, $T$ has the form
\begin{equation}%\label{eq:T}
  T=\{e\in C_0(\R_+)\mid \phi_-(t)\le e(t)\le \phi_+(t) \,\,\forall t\le \bar{t} \}\,.
 \end{equation}
By selecting a specific $\bar{e}\in X_{\underline{\theta}}$ such that $A\bar{e}=b$, and
 selecting $\phi_{\pm}(t)$ such that $\phi_{-}(t)<\bar{e}(t)<\phi_{+}(t)$ for all $t\le \bar{t}$,
 we may ensure that
\begin{equation}\label{eq:CQ:1}
b\in A( X_{\underline{\theta}}\cap  \func{int}T )
\end{equation}
which suffices for duality to hold, if the cost-function is of the form $f_0+\iota_T$
(where $f_0$ is any finite-valued functional, such as one of those listed in
Section~\ref{sectsometdperf :cont}).

Suppose also that the reference signal $w$ is minimum-phase, and that all the unstable poles and zeros
of the plant are oscillatory (that is, they all have non-vanishing imaginary part). As observed in Theorem~\ref{cor1:cont},
when $T$ is absent,  have OS${}_{ \text{\scriptsize{\funcc{opt}}} }=0$. Now introduce the further constraint represented by the set $T$,
and consider
$$
\text{OS}^T_{ \text{opt} } :=  \inf_{e\in T\cap X_{\underline{\theta}}    \cap A^{-1}b}f_{os}(e)
=  \inf_{  \overline{X_{\underline{\theta}}   } \cap A^{-1}b}(f_{os}+\iota_T )\,.
$$
(where the latter equality follows from (\ref{eq:CQ:1})).
With reference to Proposition~\ref{prop:summary}, make the additional assumption %on $\phi_{\pm}(t)$
in the case when
 $\overline{X_{\underline{\theta}}   }=C_{0,\alpha}$, that $\alpha\ge 0$, in which case we then have
 $$
 \phi_+(0)>\alpha\ge\max\{\phi_-(0) , 0 \}\,.
 $$
 (remember that in this case $\alpha=w(0+)=\bar{e}(0)$)

\begin{prop}\label{prop:tradeoff}
With the above assumptions, $\funcc{OS}^T_{ \scriptsize{\funcc{opt}} }=0$.
\end{prop}
Thus, in particular, for step-input, the imposition of rise-time constraints does not degrade the
optimal overshoot performance.

\begin{Proof}
Now, $f(e)$ is of the form $\sup_{t\ge 0}F_t(e(t))$, where
$F_t(\xi) =F^{\text{os}}(\xi)+ \iota_{\overline{(\phi_-(t),\phi_+(t))} }(\xi)$ from which follows
(since $\phi_+(t)\ge 0$) that $e_t^+(L) =\phi_+(t)$ and $e_t^-(L) =\max\{ \phi_-(t), -L\}$.
and that Assumptions 1---4 for Theorem~\ref{thmfstar1:cont} are satisfied (note, for assumption 2, needed $\phi_+(t)\ge 0$
to get $[F_t\le L]\neq\emptyset$), and hence for any $\mu\in U\oplus V\oplus W$, have
\begin{align*}
 f^*(\mu) & =  \int_0^{\bar{t}}\phi_+\,d\mu_+  +(+\infty)\cdot\mu_+((\bar{t},+\infty)) +\\
   & \qquad +\sup_{L\ge 0}\left[  \int_0^{\bar{t}}\max\{\phi_-, -L  \}d\mu_-   - L(  \mu_-((\bar{t},+\infty))+1)   \right]\\
      & \ge \int_0^{\bar{t}}\phi_+\,d\mu_+  + \int_0^{\bar{t}}\max\{\phi_-, 0 \}d\mu_-    + (+\infty)\cdot\mu_+((\bar{t},+\infty))
 \end{align*}
where the first two terms of the latter are clearly finite. Thus, whenever $\mu\in\func{dom}f^*$, it must follow that $\mu_+((\bar{t},+\infty))=0$,
and so $\func{dom}f^*\subseteq \{\mu\mid \mu|_{(\bar{t},+\infty)} \le 0\}$.

The dual then has the form (when $\overline{X_{\underline{\theta}}   }=C_{0 }$)
$$
\max_{
{\scriptsize
\begin{array}{c}
 e^*\in U\oplus V\oplus W \\
e^*\in\func{dom}f^{*} \\
e^*|_{(\bar{t},+\infty)} \le 0\\
\end{array}
} } \left[ \langle \func{Proj}_U(e^*),w\rangle -f^{*}(e^*)\right]
,
$$
but each dual $e^*$ now satisfies $e^*=0$ by Lemma~\ref{prop:sinus:cont}, by the same argument  as used in Theorem~\ref{cor1:cont}.
Hence, since $f^*(0)=\sup_{L\ge 0}-L =0$, the dual has value 0.

If, instead, $\overline{X_{\underline{\theta}}   }=C_{0,\alpha}$, the dual   takes the form
$$  (D)\,=
\max_{
{\scriptsize
\begin{array}{c}
 e^*\in U\oplus V\oplus W \\
e^*\in\func{dom}(-f^{\#}) \\
 \end{array}
} } \left[ \langle \func{Proj}_U(e^*),w\rangle +f^{\#}(e^*)\right]\,,
$$
where
$$
f^{\#}(\mu)=\max_{\lambda\in\R}(\alpha\lambda-f^*(\mu+\lambda\delta))\,.
$$
We want to show that $\func{dom}(-f^{\#})\subseteq\{ \mu  \mid \mu|_{(\bar{t},+\infty)   }\le 0   \}$.
Now, we may apply Theorem~\ref{thmfstar1:cont} for $f^*(\mu)$ for any $\mu\in U\oplus V\oplus W\oplus\R\delta$
(since for any $L\ge 0$, $t\mapsto e_t^{\pm}(L)$ has no discontinuities at $t=0$---only at $t=\bar{t}>0$, where
$|\mu|(\{ \bar{t}\})=0$). Thus,
\begin{align*}
-f^{\#}(\mu) & = \min_{\lambda\in\R}(-\alpha\lambda+f^*(\mu+\lambda\delta))\\
    & \ge \min_{\lambda\in\R} \bigg(-\alpha\lambda +\int_0^{\bar{t}} \phi_+ d(\mu+\lambda\delta)_+  \,  +\,
         \int_0^{\bar{t}}  \max\{\phi_-, 0\}d(\mu+\lambda\delta)_-  \,\,\,\,  + \\
           &  \qquad\qquad (+\infty)\cdot(\mu+\lambda\delta)_+((\bar{t},+\infty)) \bigg) \\
            & = (+\infty)\cdot\mu _+((\bar{t},+\infty)) + \int_0^{\bar{t}} \phi_+d\mu_+  \,  + \,
              \int_0^{\bar{t}} \max\{\phi_-, 0\}d\mu_-  \,\,\,\,+  \\
                 &  \qquad\qquad +\, \min_{\lambda\in\R}\big(-\alpha\lambda+ \lambda_+\phi_+(0) + \lambda_- \max\{\phi_-(0),0   \}\big)\qquad\text{ (by Lemma~\ref{lem:muFL})}\\
                 & =  (+\infty)\cdot\mu _+((\bar{t},+\infty)) + \int_0^{\bar{t}} \phi_+d\mu_+ \, + \,
              \int_0^{\bar{t}} \max\{\phi_-, 0\}d\mu_-
\end{align*}
where we used   the relation $\phi_+(0)>\alpha\ge\max\{\phi_-(0) , 0 \}$ to ensure the latter minimization has value zero.

Therefore, if $-f^{\#}(\mu)<+\infty$, then $\mu _+((\bar{t},+\infty)) =0$ so that $\mu|_{ (\bar{t},+\infty)}\le 0$, which implies,
since $\mu\in U\oplus V\oplus W$, that $\mu=0$, whence (D)$\,=f^{\#}(0)=\max_{\lambda\in\R}(\alpha\lambda-f^*( \lambda\delta))$. But, $f^*( \lambda\delta))=
\sup_{L\ge 0}(\lambda_+e_0^+(L)-L) =\lambda_+\phi_+(0)$ since $e_0^+(L)=\phi_+(0)$ for all $L$, yielding
(D)=$\max_{\lambda\in\R}(\alpha\lambda-\lambda_+\phi_+(0))=0$ as $\alpha<\phi_+(0)$ (and $\alpha\ge 0$).
Thus, $\text{OS}^T_{ \text{opt} }=0$ in this case also.
\end{Proof}

\section{Conclusions}
  Using a dual formulation, new results on
fundamental time-domain performance limitations for
continuous-time systems have been presented.   For the problem of
designing a feedback system to optimally track a specific input,
or reject a specified disturbance, there are many time-domain
performance measures of the output signal that can be used. In
addition to overshoot, undershoot and the infinity norm of the
error signal, a performance measure of practical significance,
termed fluctuation, has been investigated for the first time in a
continuous-time setting.

\bigskip

\section{Appendix A: Proof of Theorem~\ref{thmfstar1:cont}}

Firstly, we consider the case where, for each $L$, the $e_{t}^\pm(L)$ are bounded in $t\ge 0$.
Now,
$
f^*(\mu) = \sup_{L\ge 0}\sup_{e\in C_0,\,
                              f(e)\le L} (\la\mu,e\ra-L)\,,
$
noting that if $F$ is independent of $t$, the supremum is restricted to $L$ satisfying $L\ge F(0)$, as indicated in Remark~\ref{rem:thmfstar1:cont}.
We fix a value of $L$ and evaluate the inner supremum, and will show that
$$  \sup_{\scriptsize \begin{array}{c}
                                e\in C_0\\
                              f(e)\le L \end{array} } \!\!\!\la\mu,e\ra=
  \int_0^{+\infty}\!\! \!e_t^+(L)d\mu_+(t)+
 \int_0^{+\infty }\!\!\!e_t^-(L)d\mu_-(t)
$$
 from which \eqref{eq:conj:cont} then  follows.

For $\mu\in \mathbf{M}(\R_+)$, let $(P,N)$ be a Hahn decomposition
\cite{Cohn} of $\R_+$ relative to the measure $\mu$. Now, among all
$e(\cdot):\R_+\to\R$ satisfying $e(t)\in [e_t^-(L),e_t^+(L)]$ for all $t$,
$\la\mu,e\ra$ is maximized at $\bar{e}$, where
$$ \bar{e}(t):= \left\{ \begin{array}{cl}
                         e_t^+(L) & t\in P \\
                       e_t^-(L) & t\in N
                    \end{array}\right. \,.    $$
This maximal value
$\la\mu,\bar{e}\ra$ at   majorizes the required quantity.
Note that $\bar{e}$ is bounded, so $\int|\bar{e}|\,\,d|\mu|<\infty$, but generally not continuous,
nor does it decay as $t\to+\infty$. It remains to
    show that this upper bound is approached for
$e\in C_0$.

Let $\e>0$. By regularity \cite{Cohn} of the measure $|e^+-e^-|d\mu$, there are
  $K_P\subseteq P\subseteq U_P$ and
$K_N\subseteq N\subseteq U_N$ ($U$ open, $K$ compact) for which
$$
\int_{(U_P\backslash K_P) \cup   (U_N\backslash K_N)   }
\!\!\!\!\!\!\!\!\!\!\!\!\!\!\!\!\!\!\!\!\!\!\!\!\!\!|e^+-e^-|\,\,d|\mu| \le \epsilon\,.
$$
Suitable shrinkage of $U_P$, $U_N$ can ensure that $U_P\cap K_N=\emptyset=U_N\cap K_P$
(since we can take $P$ and $N$ to be disjoint).
Let $\{ \varphi_P, \varphi_N  \}$ be a partition-of-unity  (see \cite{Dug}) relative to the open covering
$\{U_P, U_N\}$ for $\R_+$. We could now propose that $\widetilde{e}(t):=e_t^+(L)\varphi_P(t)+
e_t^-(L)\varphi_N(t)$, so that $\widetilde{e}(t)\in [e_t^-, e_t^+]$ for all $t$. However,
$\widetilde{e}$ again may not be continuous. To remedy this, we seek suitable continuous
approximations for $e^+$ and $e^-$.

Let $\bar{t}>0$ be any point of (jump-)discontinuity for $e^\pm$ ($e^+$ or $e^-$ or both).
For ease of presentation, we assume strict inequality in  assumption 3.

For definiteness, assume $e_{\bar{t}-}^+> e_{\bar{t}+}^+ $. Form a line-segment, of
large negative slope $-\lambda$, from $ (\bar{t},  e_{\bar{t}+}^+)\in\R^2 $ to the point of first
 contact with the graph of $e^+|_{[0,\bar{t})  }$. Using this, we define a function $\rho^+$,
 whose {\em epigraph}\footnote{the epigraph of a function $f:X\to [-\infty, +\infty]$ is defined to be the set $\{(x,\lambda)\in X\times\R \mid f(x)\le\lambda  \}$
of points `above' the graph of $f$. }
 is the union of the epigraph of $e^+$ with the set of points above the
indicated line-segment. This yields a  function continuous at $\bar{t}$, satisfying $\rho^+\le e^+$
on $\R_+$, with $\rho^+\equiv e^+$ except on a small interval
$[ \bar{t}-\alpha(\lambda),\bar{t}]$ where $\alpha(\lambda)\to 0$ as $\lambda\to +\infty$.
Further, since $e_{\bar{t}-}^+> e_{\bar{t}+}^+ $, on taking $\lambda$ large enough, we ensure, by
the continuity of $e^-_\cdot$ near (but not at) $\bar{t}$, that also $\rho^+(t)\ge e_t^-$, for all $t$.
An analogous construction applies (using line-segments of large {\em positive\/} slope) if instead
$e_{\bar{t}-}^+< e_{\bar{t}+}^+ $.

Similarly, we obtain $\rho^-$ continuous at $\bar{t}$, for which $e^-\le  \rho^- \le e^+$ on $\R_+$
and $\rho^- \equiv e^-$ off a small interval at $\bar{t}$. It now follows that
\[
\lim_{\lambda\to +\infty }\int  |\rho^\pm - e^\pm|\,d|\mu| =
\text{constant}\cdot[\text{jump at $\bar{t}$ for }e^\pm]\cdot |\mu|(\{0\}) =0\,.
\]
Repeating this process at each such $\bar{t}$ (these are countable in number) yields continuous functions on $\R_+$,
again denoted by $\rho^\pm$, such that $e^-\le \rho^\pm \le e^+$, with
$\int  |\rho^\pm - e^\pm|\,d|\mu|\le \e/2$.

[Remark: if instead, have equality in Assumption 3.\ then the above line-segment construction will
not do---in this case one may
easily adapt this, by using nonlinear segments that are local graphs of   continuous functions.]

Now, we may define continuous functions $\widetilde{e}$, by
$ \widetilde{e}:=\rho^+ (t)\varphi_P(t)+\rho^- (t)\varphi_N(t)$. Since the $0\le \varphi_{P, N}\le 1$ and
$\varphi_P + \varphi_N \equiv 1$ we get that $e^- \le   \widetilde{e}  \le e^+$.  Also,
\begin{align*}
|\la\mu,\bar{e}\ra-\la\mu,\widetilde{e}\ra|&\le \int |\bar{e}-
\widetilde{e}|\,d|\mu| \\
& \le  \int_{(U_P\backslash K_P) \cup   (U_N\backslash K_N)   }\!\!\!\!\!\!\!\!\!\!\!\!\!
\!\!\!\!\!\!\!\!\!\!\!\!\!|e^+-e^-|\,\,d|\mu|+
 \int_{  K_P    }\!\!\!|\rho^+-e^+|\,\,d|\mu| +
\int_{ K_N   }\!\!\!|\rho^--e^-|\,\,d|\mu| \\
&\le 2\e
 \end{align*}
 At this stage, the proof is not quite complete, since
$\widetilde{e}$ may not tend to 0 as $t\to\infty$. For the final
step we form a sequence of truncations $\widetilde{e}_n$ such that:
$\widetilde{e}_n(t)=0$ for all $t\ge n+1$; they are continuous; and satisfy
$\la\mu,\widetilde{e}_n\ra \to \la\mu,\widetilde{e}\ra$ as $n\to\infty$.
We will   ensure continuity for these truncations by again using steep
interpolating line-segments.

If $\widetilde{e}(n+1)=0$, we merely take $\widetilde{e}_n(t)=0$ for $t\ge n+1$. If instead,
$\widetilde{e}(n+1)$ is nonzero, we join the point $(n+1,0)\in\R^2$ to the graph of
$\widetilde{e}|_{[0,n+1)}$ by a steep-but-nonvertical line segment, to yield a function
$\widetilde{e}_n$ that vanishes on $[n+1,+\infty)$ and agrees with $\widetilde{e}$ on $[0,n]$,
with $| \widetilde{e}_n  |\le |\widetilde{e}  |$ on $\R_+$. Also, since $0\in [F_t\le L]$ for all
$t\ge t_0$, then on taking $n\ge t_0$ and noting that $\widetilde{e}(t)\in [F_t\le L]$ for all $t$,
we have $\widetilde{e}(t)\in \func{conv}\{0, \widetilde{e}(t) \}    \subseteq[F_t\le L]$ for all $t$.
Thus, $  \widetilde{e}_n \in C_0(\R_+)$, $f(\widetilde{e}_n)\le L$, and
 \begin{align*}
\int|\widetilde{e}_n-\widetilde{e}|\,\,d|\mu| & = \int_{n }^\infty\!\!\!|\widetilde{e}|\,\,d|\mu|+
\int_{n}^{n+1} \!\!\!|\widetilde{e}_n-\widetilde{e}|\,\,d|\mu| \\
&\le 3 \int_{n}^{\infty}\!\!\! |\widetilde{e}|\,\,d|\mu|\qquad \text{ since }\,|\widetilde{e}_n |\le | \widetilde{e}|\\
& \,\,\to 0 \quad\text{ as }n\to\infty\,,
 \end{align*}
and so it follows that
$$
\la\mu,\bar{e}\ra \le \la\mu,\widetilde{e}\ra+2\e \le \la\mu,\widetilde{e}_n\ra +3\e \le
     \sup_{\scriptsize \begin{array}{c}
                                e\in C_0\\
                              f(e)\le L \end{array} }\la\mu,e\ra
+3\e\,.
$$
As $\e>0$ is arbitrary, this establishes the required equality, and hence the formula for $f^*(\mu)$, for the case
of bounded $e^\pm$.

For the unbounded case, for each $R>0$, let $F_t^R:= F_t +\iota_{[-R,+R  ]}$, where the latter denotes the indicator--function for
an interval. Then the sublevel--sets of $F_t^R$ have the form $[F_t^R\le L]=[  e_t^{R-}(L) ,  e_t^{R+}(L)    ]$ where
\[
e_t^{R+}(L) :=\min\{R, e_t^{ +}(L)\}\,, \quad\text{ and } \quad e_t^{R-}(L) :=\max\{-R, e_t^{ -}(L)\}\,.
\]
Define $f_R(e):= \sup_{t\ge 0}F_t^R(e(t))$. Then, for any allowed $\mu$ (since the  $e^{R\pm}$ are bounded and $f_R$ satisfies
Assumptions 1---4) we may apply the previous argument to $f_R$, to yield
\[
f_R^*(\mu)=\sup_{  L\ge 0  } \left(\int_0^{+\infty}\!\! \!e_t^{R+}(L)d\mu_+(t)+
 \int_0^{+\infty }\!\!\!e_t^{R-}(L)d\mu_-(t)-L\right)
\]
for any allowed $\mu$.
By an easy check, for any $e\in C_0$ and any $R>\sup|e|$,
$
f(e)=f_R(e)=\inf_{R>0} f_R (e)=\lim_{R\to\infty}f_R(e)\quad(\text{since }f_R\downarrow \text{ as } R\uparrow)
$
so that
\begin{align*}
f^*(\mu) & =\sup_{R\ge 0} f_R^*(\mu) =\sup_{L\ge 0} \,\sup_{R\ge 0}
\left(\int_0^{+\infty}\!\! \!e ^{R+}(L)d\mu_+ +
 \int_0^{+\infty }\!\!\!e ^{R-}(L)d\mu_- -L\right)\\
 & =
 \sup_{L\ge 0}
\left(\int_0^{+\infty}\!\! \!e ^{+}(L)d\mu_+ +
 \int_0^{+\infty }\!\!\!e ^{-}(L)d\mu_- -L\right)\\
\end{align*}
the latter equality (in $\R \cup\{+\infty\}$) following from the monotone convergence $e^{R+}\uparrow e^+$ and $e^{R-}\downarrow e^-$.

\section{Appendix B: Proof of Proposition~\ref{prop:summary}}

 The proof will
require a couple of preparatory lemmas. We remind the reader of
the following {\bf Notation:}
$$\R_+:= [0,\infty)$$
$$ C_0 (\R_+):=\{\varphi:\R_+\to \R\mid \varphi \text{ continuous,
} \varphi(\infty)=0 \}$$
$$ C_{0,\alpha}(\R_+):=\{\varphi\in C_0(\R_+)\mid\varphi(0)=\alpha \}$$
%$$ C_b (\R_+):=\{\varphi \text{ continuous and bounded on }
%\R_+ \}$$
$$
\underline{\theta}:= (\theta_p, \theta_w),\quad\text{where}\quad \theta_p\geq
0,\,\,\theta_w\geq 1; \quad\text{ and }\quad
\theta:=\theta_p+\theta_w
$$

 It will
 be seen that the closure of $X_{\underline{\theta}}$ is always one of
$C_0$ or $C_{0,\alpha}$ (for suitable $\alpha$), depending on the
value of $\theta_w,\theta_p$.

\bigskip
%\subsection{Case of $\theta_w\geq 1$}

We start with   a simple result for Laplace integrals.

\begin{lem}\label{lem:A1} Let $\varphi:\R_+ \to \R$ such that
  $\varphi^{(k)}(t)=O(e^{at})$ (for some $a>0$)
and $\varphi^{(i)}(0+)$ exists and finite, for $i=0,1,\ldots k$.
Then for all $s\in \C$ for which $\func{Re}s>a$,
$$
\widehat{\varphi^{(k)}}(s)=\widehat{\varphi}(s)s^k-\varphi(0+)s^{k-1}
-\varphi^{'}(0+)s^{k-2} - \ldots - \varphi^{(k-2)}(0+)s - \varphi^{(k-1)}(0+) \,.
%\mbox{ for } \mbox{Re }s>a\,.
$$  \end{lem}
\begin{Proof}   Integrating by
parts yields
$\widehat{\varphi^{(k)}}(s)=s\widehat{\varphi^{(k-1)}}(s)-
\varphi^{(k-1)}  (0+)$ for such $s$. Repeat for $\varphi^{(k-1)}$
etc.
\end{Proof}

\smallskip

We   assume that the function $w$   has rational Laplace
transform. This ensures that $\widehat{w^{(i)}}(s)\to 0 $ as
$\text{Re }s\to +\infty$, and $w^{(i)}(0+)=w^{(i)}(0)$, for all
$i=0,1,2 \ldots $. Moreover, Lemma~\ref{lem:A1} applies, for every
$k$.

Now, define
\begin{align*}
X_{\underline{\theta}} :=\{ e\in C_0(\R_+)\mid
\hat{e} & \text{ rational, } \hat{e}(\infty)=0 \text{ (order at least }\theta_w) \\
&  \widehat{(e-w)}(\infty)=0 \text{ (order at least
}\theta=\theta_w +\theta_p) \}.
\end{align*}

Also, form the subspace
$$
\begin{array}
[c]{cc}%
Y_{\underline{\theta}}:= & \left\{ e:  t\mapsto\text{Re}\sum_{i\in I}c_i
t^{k_i}e^{-\lambda_i t}\,\bigg{|}\,
\begin{array}
[l]{l}%
    c_i\in \C,\,k_i\geq \theta_w-1,\, \text{Re
}\lambda_i>0,\,
  I  \text{ finite,}   \\
    e^{(j)}(0)=w^{(j)}(0+), \,j=0,1,\ldots ,\theta-2
\end{array}
\right\}
\end{array}
$$

\begin{lem}\label{lem:A2}\quad $Y_{\underline{\theta}}\subseteq
X_{\underline{\theta}}$.
\end{lem}
\begin{Proof} If $e\in Y_{\underline{\theta}}$, then clearly
$\hat{e}$ has a zero of order $\theta_w$ at infinity. Now,
$[(e-w)^{(\theta-1)}]\,\widehat{}\, (s)\to 0$ when Re
$s\to\infty$. By Lemma~\ref{lem:A1} applied to $\varphi:=e-w$, we
obtain, when $s\to\infty$,
\begin{align*}
 0\leftarrow [(e-w)^{(\theta-1)}]\,\widehat{}\, (s)
&=s^{\theta-1}(\hat{e}(s)-\hat{w}(s))-s^{\theta-2}(e-w)(0+)-\ldots
-(e-w)^{(\theta-2)}(0+)\\ & =s^{\theta-1}(\hat{e}(s)-\hat{w}(s))
\end{align*}
 on
using the constraints at $t=0+$. Hence $\hat{e}-\hat{w}$ has a
$\theta$--order zero at $\infty$.
\end{Proof}

\begin{lem}\label{lem:A3} If
 $e\in X_{\underline{\theta}}$, then:
$e^{(j)}(0)=0$ for $j=0,\ldots ,\theta_w\!-2$; and also
$e^{(j)}(0)=w^{(j)}(0+)$ for $j=0,\ldots ,\theta-2$.
\end{lem}
\begin{Proof}
Let $k:=\theta\!-\!2$. Now $[(e-w)^{(k+1)}]\,\widehat{}\, (s)\to
0$ at infinity. From Lemma~\ref{lem:A1},
$s^{k+1}(\hat{e}(s)-\hat{w}(s))-(e(0)-w(0+))s^k-\ldots -
(e^{(k)}(0)-w^{(k)}(0+))=[(e-w)^{(k+1)}]\,\widehat{}\, (s)\to 0$.
By assumption on $e$, have $s^{k+1}(\hat{e}-\hat{w})(s)\to 0$ at
infinity, and this forces the equalities $e(0)=w(0+),\ldots
,e^{(k)}(0)=w^{(k)}(0+)$. For $j=0,\ldots, \theta_w \! -\!2$,
argue as above using $e$ instead of $e-w$, to get $e^{(j)}(0)=0$.
\end{Proof}

\begin{remark}\label{rem:1} If $\theta_w\ge 2$, then Lemma~\ref{lem:A3}
requires $w^{(j)}(0+)=0$ for $j=0,\ldots,\theta_w-2$, otherwise
the statement  of the lemma entails a contradiction, and
$X_{\underline{\theta}}$ will be empty, hence of no interest. If
$\theta_w=1$, the   conditions $e^{(j)}(0)=0$ ($j\le \theta_w-2$)
are absent, so $w(0+)$ is unrestricted in this case.
\end{remark}

 For any $\underline{\theta}$, form the set
 $W_{\underline{\theta}}\supseteq Y_{\underline{\theta}}$ by
 removing the constraints at $t=0$ in the latter, that is,
$$
W_{\underline{\theta}}:=\left\{e(t)=\text{Re}\sum_{i\in I}c_i
t^{k_i}e^{-\lambda_i t}\, \bigg{|}\,  c_i\in \C,\,k_i\geq
\theta_w-1,\, \text{Re }\lambda_i>0,\,
  I  \text{ finite}
   \right\}.
$$

Observe that $e(t)=O(t^m)$ at $t=0$ (where $m:=\theta_w-1$), it
follows that $e^{(j)}(0)=0$ for $j=0,1,\ldots,m-1$ when
$\theta_w>1$.

\bigskip

We turn now to the proof of the Proposition. Because of the
preceding lemmas, it suffices to prove density for the subsets
$Y_{\underline{\theta}}$ of $X_{\underline{\theta}}$.

\subsubsection{Case of $\theta_w> 1$}

If $\theta_w>1$, note that $\varphi(0)=0=\varphi(+\infty)$
whenever $\varphi\in W_{\underline{\theta}}$. Also,
$W_{\underline{\theta}} $ forms a  algebra of continuous functions
on the {\em one-point compactification\/} (see \cite{Dug})
$[0,+\infty]$ of $\R_+$, since
$\varphi(+\infty)=\lim_{t\to\infty}\varphi(t)$ exists for each
such $\varphi$. We seek to apply the Stone--Weierstrass Theorem
(see, for instance, \cite{Dug}, or any functional analysis text)
to establish the density of $W_{\underline{\theta}}$ in a suitable
subspace of $C([0,+\infty])$. However, there is a difficulty, in
that $W_{\underline{\theta}}$ does not `separate all points' of
$[0,+\infty]$ in the sense of the Theorem---indeed, by definition
of $W_{\underline{\theta}}$,  the points $0$ and $+\infty$ cannot
be so separated (if $\varphi\in W_{\underline{\theta}}$, it
attains equal values at $0$ and at $+\infty$). This is the only
pair of points that cannot be separated in this sense. We resolve
this issue by insisting that $0$ and $+\infty$ are `the same
point'.

To do this, let $X$ denote the   the quotient   space, formed by
identifying the points $0$ and $+\infty$ in $[0,+\infty]$, and
endowing this with the quotient topology \cite{Dug}. Then $X$ is
also compact, being the image of $[0,+\infty]$ under the quotient
map: $[0,+\infty]\to X$, which is continuous by definition of the
quotient topology. It is easily verified that $X$ is also
Hausdorff. Each member of $W_{\underline{\theta}} $ may now be
viewed as a continuous function on $X$, since it has a
well-defined value (zero) at $[0]=[+\infty]\in X$.

Clearly, $W_{\underline{\theta}}\oplus\R :=\{\varphi + c\!\mid \!\varphi\in W_{\underline{\theta}},\,\,
c\in\R\} $ forms a subalgebra of $C(X)$ that contains the constant
functions, and now separates the points of $X$ (since the
offending pair in $[0,+\infty]$ have been merged). We may now
apply Stone--Weierstrass, to conclude that
$W_{\underline{\theta}}\oplus\R$ is dense in $C(X)$ under the
supremum--norm.

Returning to $[0,+\infty]$, this means that any $\varphi\in
C([0,+\infty])$ for which $\varphi(0)= \varphi(+\infty)$ may be
uniformly approximated by functions from
$W_{\underline{\theta}}\oplus\R$, whence follows that
$\overline{W_{\underline{\theta}} }=C_{0,0}(\R_+)$.

\smallskip

\smallskip

To finally prove the required property, let $\varphi\in
C_{0,0}(\R_+)$, and let $\e
>0$. By the   density of  $W_{\underline{\theta}}$, there is $\widetilde{\varphi}\in
W_{\underline{\theta}}$ such that
$\|\widetilde{\varphi}-\varphi\|_\infty \le \e $. If $\theta_p=0$,
we are done, since in this case,
$W_{\underline{\theta}}=Y_{\underline{\theta}}$. If $\theta_p>0$,
then   $\widetilde{\varphi}$ may not satisfy the conditions at
$t=0$ required for membership of the subset
$Y_{\underline{\theta}}$, so we need to perturb
$\widetilde{\varphi}$ to achieve this. Place
$$
\xi :=  \left( w^{  }(0+)-\widetilde{\varphi}^{  }(0),\ldots,
w^{(\theta_w+\theta_p-2)
}(0+)-\widetilde{\varphi}^{(\theta_w+\theta_p-2 ) }(0)\right)
=(0,\ldots,0,\xi_{\theta_w-1}, \ldots,\xi_{\theta_w+\theta_p-2  })
$$
the latter obtaining since $\widetilde{\varphi}^{(j ) }(0)=0$ for
$j\le \theta_w-2$ and also if $X_\theta\neq\emptyset$ we must have
$w^{(j)  }(0+)=0$ for such $j$.

We now seek $p(\cdot)\in W_{\underline{\theta}}$ such that
$\widetilde{\varphi}+p\in Y_{\underline{\theta}}$ and
$\|p\|_\infty$ is suitably small. We slightly modify a procedure
used
 in \cite{Miller}. Try $p(t)=t^m p_0(t)$ where $m=\theta_w-1$, and
$p_0(t)=\sum_{l=0}^q c_l e^{-(l+1)\sigma t}$ with $q:=\theta_p-1$,
and $\sigma >0$  so large that $\sup_{t\ge 0}t^m e^{-\sigma t}\le
\e/(1+|\xi|)$. Then, for $k$ such that $m\le k\le
m+q=\theta_w+\theta_p-2$, have $p^{(k)}(0)=
[k!/(k-m)!]p_0^{(k-m)}(0)$ (as well as $p^{(k)}(0)=0$ for $k< m$).
For $\widetilde{\varphi}+p$ to satisfy the required constraint at
$t=0$, we need
$$p_0^{(k-m)}(0)=[(k-m)!/k!]\xi_k\qquad k=m,m+1,\ldots,m+q$$
or, equivalently,
$$
\left[\begin{array}{cccc}
1 & 1 & \dots & 1\\
1 & 2 & \dots & q+1\\
1 & 2^2 & \dots & (q+1)^2\\
\vdots & \vdots & \ddots & \vdots\\
1 & 2^q & \dots  & (q+1)^q
\end{array}\right]
\left[\begin{array}{c}
c_0 \\
c_1 \\
\vdots \\
c_q
\end{array}\right]
=  \left[\begin{array}{c}
 0!\xi_m/m! \\
1! \xi_{m+1}/(-  (m+1)!\sigma)\\
 \vdots\\
 q!\xi_{m+q}/((-1)^q  (m+q)!\sigma^q)
\end{array}\right]
$$
Let $\gamma$ denote a suitable norm of the inverse of the above
square matrix. Place $K(q)=(\gamma(q)+1)(q+1)$---note that it
depends {\em only\/} on  $q$. Solving the above for the
coefficients $c_l$ yields $p_0$ satisfying $ |p_0(t)|\le
K(q)|\xi|e^{-\sigma t} $ for all $t\ge 0$, so that
$$
|p(t)|\le K(q)t^m e^{-\sigma t}|\xi|\le K(q)\e \quad\text{ for all
}t\ge 0\,.
$$
Thus it follows for this choice of $p$ that
$\widetilde{\varphi}+p\in Y_{\underline{\theta}}$, and
$$
\|\widetilde{\varphi}+p-\varphi\|_\infty \le \e + \|p\|_\infty \le
(1+K(q))\e\,.
$$
Since $K$ does not depend on the choice of $\varphi$ or
$\widetilde{\varphi}$ etc, it follows from the arbitrariness of
$\e>0$ that $Y_{\underline{\theta}}$, and hence
$X_{\underline{\theta}}$, is dense in $C_{0,0}(\R_+)$, as claimed.

\subsubsection{Case of $\theta_w= 1$}
The same reasoning applies, but without recourse to quotient
spaces of $[0,+\infty]$.

When $p>0$, we obtain that $Y_{\underline{\theta}}$ is dense in
$C_{0,\alpha}(\R_+)$, with $\alpha=w(0+)$.

 If $p=0$, the constraints at $t=0$ in the definitions of
$Y_{\underline{\theta}}$ and $W_{\underline{\theta}}$ are absent,
and arguing as above shows that $Y_{\underline{\theta}}$ is dense
in the entire space $C_0(\R_+)$.

\bigskip

This completes the proof of Proposition~\ref{prop:summary}.

\section {Appendix C:  Proof of Lemma~\ref{lem:fhash}}

 We start with
a lemma.
\begin{lem}\label{lem:mu}
Let $\mu\in {\bf M}(\R_+)$. Then, $$ \|\mu
+\lambda\delta\|=\|\mu\|+|\mu(\{0\})+\lambda|-|\mu(\{0\})|\,.$$ In
particular, if $\mu(\{0\})=0$ (for example, when $\mu\in U\oplus
V\oplus W$) we have $ \|\mu +\lambda\delta\|=\|\mu\|+ |\lambda| $.
\end{lem}
\begin{Proof} By definition,
\begin{align*}
\|\mu +\lambda\delta\|&=|\mu +\lambda\delta|(\R_+)=\sup
\left\{\sum_{I\in\Pi}|\mu +\lambda\delta(I)| : \Pi \text{
partitions $\R_+$ into intervals} \right\}\\
&=\lim_{\|\Pi\|\to 0}\left( \sum_{I:0\notin I}|\mu(I)|  +
|\mu(I_0)+\lambda|\right)
\end{align*}
where $I_0$ denotes the unique member of $\Pi$ that contains $0$.
Since $I_0\downarrow \{0\}$ for a subsequence of partitions, we
have $\mu(I_0)\to \mu(\{0\})$, so that
\begin{align*}
\|\mu +\lambda\delta\|&= \lim_{\|\Pi\|\to 0}\left( \sum_{I \in
\Pi}|\mu(I)|  + |\mu(I_0)+\lambda|-|\mu(I_0)|\right)\\
&=\|\mu\|+|\mu(\{0\})+\lambda|-|\mu(\{0\})|
\end{align*}
\end{Proof}

Now, consider $f=f_{os}$. Then,
\begin{align*}
f_{os}^{\#}(\mu) &= \max_ {\lambda\in{\scriptsize
\R}}\,[\alpha\lambda-f_{os}^*(\mu+\lambda\delta)]\\
&=\max\{\alpha\lambda\mid \mu+\lambda\delta\le
0,\,\,\|\mu+\lambda\delta\|\le 1\}\\
&=\alpha\max\{\lambda \mid \mu+\lambda\delta\le 0,\,\,\|\mu\|+
|\lambda|\le 1\} \text{ by Lemma~\ref{lem:mu}}
\end{align*}
Let $\mu\in U\oplus V\oplus W$.  If $\|\mu\|>1$, then $\|\mu
+\lambda\delta\|=\|\mu\|+ |\lambda|>1$ for each $\lambda$, so
$f_{os}^*(\mu+\lambda\delta)=+\infty$ for all $\lambda$, and
$f_{os}^{\#}(\mu)=-\infty=-f_{os}^*(\mu)$. If $\mu$ is not a
negative measure, then $\mu(E)>0$ for a Borel set $E$ in $\R$
(with $0\notin E$), implying for all $\lambda\in\R$, that
$(\mu+\lambda\delta)(E)=\mu(E)>0$, so $\mu+\lambda\delta$ also not
negative, and again, $f_{os}^*(\mu+\lambda\delta)=+\infty$ for
such $\lambda$, yielding
$f_{os}^{\#}(\mu)=-\infty=-f_{os}^*(\mu)$.

If $\mu\le 0$ and $\|\mu\|\le 1$ (so $\mu\in\func{dom}f_{os}^*$),
then, since $\mu+\lambda\delta\le 0$ implies $\lambda\le 0$ (take
$E=\{0\}$), it follows that $f_{os}^{\#}(\mu)=0$, and hence equals
$-f_{os}^{*}(\mu)$.

This proves that $f_{os}^{\#} =-f_{os}^*$ on $ U\oplus V\oplus W$.
The proof that $f_{us}^{\#} =-f_{us}^*$ is similar.

\smallskip

To treat $f_{fl}$, we note an elementary lemma.
\begin{lem}\label{lem:muFL}
Let $\mu\in {\bf M}(\R_+)$ with $\mu(\{0\})=0$. Then,
$(\mu+\lambda\delta)_\pm= \mu_\pm + \lambda_\pm\delta$ for any
$\lambda\in\R$.
\end{lem}

Then, by similar arguments to those above,
$f^{\#}_{\text{fl}}(\mu)=-\infty$ if $\mu\notin\func{dom}f^{*}_{\text{fl}}=
\{\mu\mid \mu_+(\R_+)\le 1/2 \, \, \& \,\,\mu_-(\R_+)\ge -1/2\}$,
and for $\mu\in\func{dom}f^{*}_{\text{fl}}$,
$f^{\#}_{\text{fl}}(\mu)=\alpha(1/2-\mu_+(\R_+))$. Thus the required
relation for $f_{\text{fl}}$ is also established.

\smallskip

Similar reasoning applies for the remaining functionals.

\end{document}